\documentclass{amsart}
\usepackage[utf8]{inputenc}
\usepackage{amsmath}
\usepackage{amsfonts}
\usepackage{amssymb}
\usepackage{amsthm}
\usepackage{graphicx}
\usepackage[colorinlistoftodos]{todonotes}
\usepackage[colorlinks=true, allcolors=blue]{hyperref}
\usepackage{enumerate}
\usepackage{bm}
\usepackage{enumitem}

\newtheorem{thm}{Theorem}[section]
\newtheorem{lem}[thm]{Lemma}
\newtheorem{cor}[thm]{Corollary}

\newtheorem{prop}[thm]{Proposition}
\newtheorem{rmk}[thm]{Remark}

\theoremstyle{definition}

\newtheorem{claim}[thm]{Claim}

\title[Morse Index Bound and Morse inequalities]{Morse Index Bound of simple closed geodesics on $2$-spheres and strong Morse Inequalities}
\author{Dongyeong Ko}
\address{Department of Mathematics, Rutgers University - New Brunswick, Piscataway, NJ 08854}
\email{dk954@math.rutgers.edu}
\begin{document}
\maketitle
\begin{abstract}
 We give a Morse-theoretic characterization of simple closed geodesics on Riemannian $2$-spheres. On any Riemannian $2$-sphere endowed with a generic metric, we show there exists a simple closed geodesic with Morse index $1$, $2$ and $3$. In particular, for an orientable Riemannian surface we prove strong Morse inequalities for the length functional applied to the space of simple closed curves.

\end{abstract}
\section{Introduction}
The celebrated theorem of Lusternik-Schnirelmann \cite{LS} is:

\begin{thm}
A closed Riemannian $2$-sphere $(S^{2},g)$ admits at least three simple closed geodesics.
\end{thm}

Birkhoff \cite{B} initiated the construction of a closed geodesic by introducing the min-max method to find unstable geodesics, in case of where minimization techniques do not work. Later, Lusternik and Schnirelmann \cite{LS} produced three closed geodesics by using higher parameter families of simple closed curves from the nontrivial homology classes of the space of embedded curves. The existence of three simple geodesics comes from the topology of the space of simple curves, which can be deformed into $(\mathbb{R}P^{3} \setminus D^{3}, \partial)$. It was not clear that whether the geodesics are embedded in Lusternik and Schnirelmann's work. The proof was repaired by Ballmann \cite{Ba} and Grayson \cite{Gr} later. In particular, Grayson \cite{Gr} used curve shortening flow to prove the existence of simple closed geodesics. On the other hand, Pitts \cite{P1} developed min-max method and produced stationary geodesic networks with possible point singularities, but it does not give full regularity of critical geodesics.

In higher dimension, for $3 \le n+1 \le 7$, Morse index bounds of minimal hypersurfaces $\Sigma^{n}$ on closed manifolds $M^{n+1}$ were obtained by Marques and Neves. Marques and Neves \cite{MN2} proved that the Morse index of minimal hypersurface is bounded above by the number of parameters of families of cycles. They also settled that generically there is a lower bound of Morse index by the number of parameters under the mulitplicity one assumption in \cite{MN3}. By combining these with Zhou's Multiplicity One Theorem \cite{Z}, for bumpy metrics, there are minimal hypersurfaces with Morse index $k$ for each $k \in \mathbb{N}$. Recently, Marques, Montezuma and Neves \cite{MMN} proved the strong Morse inequalities for the area functional in codimension $1$.

 A natural question is whether there is any Morse-theoretic characterization of simple closed geodesics, while the Morse theory of closed immersed geodesics has been developed in many literature (For instance, see Chapter VII of Morse \cite{M}). From Lustenik-Schnirelmann category theory, claiming that the number of critical point of a smooth real-valued function defined on a manifold $M$ is bounded from below by one plus the maximal cup-length of the cohomology ring of $M$, it is expected that there are three simple closed geodesics with index $1,2,3$ from each parameter family on a bumpy sphere. Moreover, it has been known that there are exactly three simple closed geodesics with Morse index $1,2$ and $3$ on ellipsoids in Morse's example (Theorem 2.1 of the Chapter IX of \cite{M}).

However, to the author's knowledge, the index estimate was unknown for simple closed geodesics on generic $2$-spheres. The main difficulty comes from the lack of the weak convergence such as the Hilbert manifold structure of the space of simple closed curves. Moreover, it is even difficult to construct a smooth homotopy between two smooth families of simple closed curves with a controlled length. While the existence of such a homotopy without length bound is known by the work of Smale \cite{Sm}, there was no such construction of homotopy with quantitative bounds. Note that Ketover, Liokumovich and Song \cite{KLS} constructed the smooth interpolation with controlled area between two surfaces which are close to a union of strictly stable minimal surfaces.

In this paper, we confirm the heuristic above by deriving Morse index bound of simple closed geodesics of generic $2$-spheres and strong Morse inequalities of closed geodesics without self-intersection for length functional. In particular, we realize the geodesics with Morse index $1,2,3$ on generic Riemannian $2$-spheres in Theorem 1.2 below.

A Riemannian metric $g$ is called \emph{bumpy} if every closed geodesic is nondegenerate i.e. there is no closed geodesic that admits a non-trivial Jacobi field. The bumpy condition gives Morse property to the length functional. Abraham showed that bumpy metrics are generic in the $C^{r}$-Baire sense in Abraham \cite{Ab} for $r \ge 5$.

\begin{thm}
For a $2$-sphere with a bumpy metric $(S^{2},g)$, for each $k=1,2,3$, there exists a closed and embedded geodesic $\gamma_{k}$ with
\begin{equation}
    index(\gamma_{k}) = k
\end{equation}
and the lengths $\gamma_{1}, \gamma_{2}$ and $\gamma_{3}$ satisfy
\begin{equation*}
    |\gamma_{1}|<|\gamma_{2}| <|\gamma_{3}|.
\end{equation*}
\end{thm}
\begin{rmk}
A triaxial ellipsoid whose lengths of three axes are slightly distinct has only three simple closed geodesics with Morse index $1,2$ and $3$ (Theorem 2.1 of the Chapter IX in Morse \cite{M}).
\end{rmk}

For general $2$-spheres, we have a following Morse-theoretic characterization as a direct corollary of Theorem 1.2:
\begin{cor}
On a Riemannian $2$-sphere $(S^{2},g)$, for each $k=1,2,3$, there exists a closed and embedded geodesic $\gamma_{k}$ with
\begin{equation*}
    index(\gamma_{k}) \le k \le index(\gamma_{k}) + nullity (\gamma_{k}).
\end{equation*}
\end{cor}
\begin{rmk}
    After this paper was completed, Stephan Suhr brought to the author's attention that Theorem 1.2 and Corollary 1.4 may be obtained from combining Theorem 1.3(iii) in De Philippis-Marini-Mazzucchelli-Suhr \cite{DMMS} and the classical arguments on infinite-dimensional Morse theory, as discussed in \cite{Ch} (See also \cite{M}, \cite{AM} and \cite{BL}).
\end{rmk}
We generalize the Morse index bound (1) by obtaining Morse inequalities for length functional for all simple closed geodesics on orientable surfaces $(M^{2},g)$ endowed with bumpy metric, even for those which do not realize the width of three nontrivial homologies of the space of curves. We also note that the example of Colding and Hingston \cite{CH} that contains simple closed geodesics with arbitrary Morse Index on orientable surfaces.

We define $b_{k}(a, \Pi)$ as the $k$-th Betti number of the space of simple closed curves whose length is less than $a$ with coefficients in $\mathbb{Z}_{2}$ for a fixed homotopy class $\Pi$ of curves, and $c_{k}(a,\Pi)$ as the number of simple closed geodesics on $(M^{2},g)$ with length less than $a$ within $\Pi$. Also let us count point curves as one geodesic with zero length and index $0$.

\begin{thm} 
For each $a \in (0,\infty)$, $b_{k}(a,\Pi) < \infty$ for every $k \in \mathbb{Z}_{+}$ and the Strong Morse inequality for simple closed geodesics in a fixed homotopy class $\Pi$ hold:
\begin{equation*}
    c_{r}(a, \Pi) - c_{r-1}(a, \Pi) +\hdots+ (-1)^{r}c_{0}(a,\Pi) \ge b_{r}(a,\Pi) - b_{r-1}(a,\Pi) +...+ (-1)^{r}b_{0}(a,\Pi) 
\end{equation*}
for every $r \in \mathbb{Z}_{+}$. In particular,
\begin{equation*}
    c_{r}(a,\Pi) \ge b_{r}(a,\Pi)
\end{equation*}
for every $r \in \mathbb{Z}_{+}$.
\end{thm}

By compactness theorem, we are able to deduce that the number of simple closed geodesics with bounded length is finite and so we can sum up all terms and we have the following:

\begin{cor} For each $a \in (0,\infty)$, $b_{k}(a) < \infty$ for every $k \in \mathbb{Z}_{+}$ and the Strong Morse inequality for simple closed geodesics hold:
\begin{equation*}
    c_{r}(a) - c_{r-1}(a) +\hdots+ (-1)^{r}c_{0}(a) \ge b_{r}(a) - b_{r-1}(a) +...+ (-1)^{r}b_{0}(a) 
\end{equation*}
for every $r \in \mathbb{Z}_{+}$. In particular,
\begin{equation*}
    c_{r}(a) \ge b_{r}(a)
\end{equation*}
for every $r \in \mathbb{Z}_{+}$.
\end{cor}

The novel idea of our proof is to construct a smooth interpolation between two families of curves near a fixed geodesic $\gamma$ satisfying a length constraint. More specifically, for small $\epsilon>0$ if $\Phi : X \rightarrow \mathcal{S}$ is a continuous map in the smooth topology such that
\begin{equation*}
    \sup \{ F(\Phi(x), \gamma) : x \in X \} < \epsilon,
\end{equation*}
then we can construct a homotopy $H:[0,1] \times X \rightarrow \mathcal{S}$ between $\Phi(x)$ and $\gamma$ such that the following $F$-distance estimate holds along $H$:
\begin{equation}
    \sup \{ F(H(t,x),\gamma) : x \in X \text{ and } t \in [0,1] \} < C(|\gamma|) \sqrt{\epsilon},
\end{equation}
where $\mathcal{S}$ is a space of simple closed curves on $(S^{2},g)$ and for some $C(|\gamma|)>0$. By deforming $H$ with a local min-max deformation originated from White \cite{W2} and taking a pushforward homotopy, we prove the interpolation lemma (Lemma 6.4) which is a smooth analog of Theorem 3.8 in \cite{MN3}. This follows the length bound along the homotopy which do not exceed the widths in the proof of Theorem 1.2.

Our construction of the homotopy $H$ in (2) comes from the analysis of homotopies in the squeezing lemma (Lemma 5.4). For small $h$ and $\epsilon$, one obtains the upper bound of $F$-distance between a simple closed curve $\alpha$ lying on $N_{h}(\gamma)$ whose length is bounded by $|\alpha|<|\gamma|+\epsilon$ and a strictly stable geodesic $\gamma$ with a negative Gaussian curvature on $\gamma$. We see this by obtaining the bound of total angle $\int_{\alpha} |\theta|$ in terms of $h$ and $\epsilon$ inspired by the proof of Quantitative Constancy Theorem in Song and Zhou \cite{SZ}. This gives rise to the idea to control the $F$-distance with $\gamma$ along the homotopy even if curves are not graphical.

Since the map constructed in Lemma 5.4 comes out from the curve shortening flow and the squeezing map in this case, the curves along the homotopy satisfy the length bound and stay in the tubular neighborhood by the existence of a mean convex foliation by level curves of the geodesic. This gives the $F$-distance bound along the homotopy. For the general geodesic case, one can obtain the upper bound of $F$-distance by considering appropriate deformation to change the geodesic to be strictly stable and with negative ambient curvature case as in the strictly stable and with negative ambient curvature case in the tubular neighborhood which preserves the order of the $F$-distance bound along the homotopy.

We adopt smooth min-max construction of geodesics as in Grayson's work \cite{G} on curve shortening flow and develop Morse-theoretic characterization of simple closed geodesics arising from smooth sweepouts with min-max techniques (See Colding-De Lellis \cite{CD} and Haslhofer-Ketover \cite{HK} for $3$-dimensional manifolds). We prove the Morse index bound by relying on the analysis of local min-max structure near the geodesic with given Morse index of White \cite{W2}, deformation techniques in Marques and Neves' work on Morse theory for area functional in \cite{MN2} and \cite{MN3} and the interpolation lemma we explained above. Moreover, for the proof of Morse inequalities, we mostly follow the schemes of Marques, Montezuma and Neves \cite{MMN} and we apply the interpolation lemmas (Lemma 6.3 and Lemma 6.4) in place of the interpolation lemma in the flat topology (Theorem 3.8 in \cite{MMN}).

The organization of this paper is as follows. In Section 2, we introduce the smooth min-max construction of geodesics and prove the pull-tight properties. In Section 3, we prove that the relation between $F$-distance and Hausdorff distance. In Section 4, we show the $1$-varifold version of local min-max theorem. In Section 5, we describe the perturbation of the metric and prove squeezing lemma on a family of curves.  In Section 6, we prove $F$-distance estimate along the squeezing map and the interpolation lemma. In Section 7, we prove Theorem 1.2. In Section 8, we obtain strong Morse inequalities for simple closed geodesics for length functional.

In Appendix, we prove the compactness theorem of geodesics which is $1$-dimensional analog of Sharp \cite{Sh}.
\section*{Acknowledgments}
The author would like to thank his advisor Prof. Daniel Ketover for suggesting this problem, valuable discussions, careful comments on the previous versions of this paper, and his constant support.  The author is thankful to Prof. Yevgeny Liokumovich for his interest in this work and introduction to related contexts and questions. The author thanks to Prof. Stephan Suhr for pointing out the paper \cite{DMMS} and detailed explanations on related works.  The author was partially supported by NSF grant DMS-1906385.

\section{Smooth Min-max Construction}

We introduce the smooth min-max setting and generate simple closed geodesics. We will construct the deformation to prove Theorem 1.2 in Section 7 based on the tightened sequence in Theorem 2.3. The following contents are smooth min-max constructions of simple closed geodesics (see \cite{CD} and \cite{HK} for minimal surface version) in myriad of classical literature (e.g. \cite{Gr} and \cite{LS}). 

First let $V_{1}(S^{2})$ be a space of varifolds and $IV_{1}(S^{2})$ be a space of integral varifolds on $(S^{2},g)$ and endow an $F$-metric on the space of varifolds.

We recall the following definition of varifold $F$-distance $F: V_{1}(S^{2}) \times V_{1}(S^{2}) \rightarrow \mathbb{R}$ in the space of varifolds on $(S^{2},g)$ from \cite{P2}.
\begin{equation}
    F(V,W) = \sup \{ V(f) - W(f) : f \in C_{c}(G_{1}(S^{2})), where\,\,|f| \le 1 , Lip(f) \le 1 \},
\end{equation}
where $G_{1}(S^{2})$ is a Grassmannian on the sphere $(S^{2},g)$ and the metric on the Grassmannian bundle $G_{1}(S^{2})$ is defined by the induced metric by $g$.

We define the space of smooth embedded closed (unparametrized) curves $\mathcal{X}$ on $(S^{2},g)$ as following:
\begin{equation*}
    \mathcal{X} := \{ F(S^{1}) \,\, | \,\, F: S^{1} \rightarrow S^{2} \text{ is a smooth embedding} \}.
\end{equation*}
Also we define the space of boundary point curves as $\mathcal{X}_{0}$. Lusternik-Schnirelmann theory (Appendix A.3 of \cite{K}) shows that the space of embedded circles retracts onto the space of geodesic circles of sphere, so we have $(\mathcal{X} \cup \mathcal{X}_{0}, \mathcal{X}_{0})$ retracts onto $(\mathbb{R}P^{3} \setminus D^{3}, \partial)$, where $D^{3}$ is a open $3$-ball. Let us endow smooth topology in $\Sigma:= \mathcal{X} \cup \mathcal{X}_{0}$ and identify boundary point curves in $\mathcal{X}_{0}$. Then since $\Sigma$ is homeomorphic to $\mathbb{R}P^{3} \setminus D^{3}$, $\mathcal{S}:= \Sigma / \partial \Sigma$ is homeomorphic to $\mathbb{R}P^{3}$. 

From the topological structure of $\mathcal{S}$ above, it has the following three nontrivial $\mathbb{Z}_{2}$-homology classes $\{ h_{i} \}_{i=1,2,3}$ (cf. Section 2 of \cite{HK}) :
\begin{equation*}
    h_{i} := H_{i}(\mathcal{S}, \mathbb{Z}_{2}) = \mathbb{Z}_{2}.
\end{equation*}
Now we define sweepouts more precisely. Let $\alpha$ be a generator of the first cohomology ring $H^{1}(\mathcal{S}, \mathbb{Z}_{2})$. We consider the cohomology ring, which is given by $H^{*}(\mathcal{S},\mathbb{Z}_{2}) = \mathbb{Z}[\alpha]/\alpha^{4}$. For each $i$, denote $X$ to be an $i$-dimensional simplicial complex. We say that $\Phi:X \rightarrow \mathcal{S}$ detects $\omega \in H^{i}(\mathcal{S}, \mathbb{Z}_{2})$ if
\begin{equation*}
    \Phi^{*} (\omega) \neq 0.
\end{equation*}
Then we let $\Phi$ be an \emph{$i$-sweepout} with an endowed smooth topology if it detects $i$-cup product $\alpha^{i}$.  Now let $S_{i}$ be the set of all $i$-sweepouts which detect $\alpha^{i}$ for $i \in \{ 1,2,3 \}$. We define the \emph{width} of $i$-parameter sweepouts as
\begin{equation*}
    \omega_{i}(S^{2}) : = \inf _{\Phi \in S_{i}} \sup_{x \in X} |\Phi(x)| = L_{i},
\end{equation*}
for $i \in \{ 1, 2, 3 \}$. By definition of $i$-sweepouts, we have $\omega_{1}(S^{2}) \le \omega_{2}(S^{2}) \le \omega_{3}(S^{2})$.

Let us denote the sequence of $i$-sweepout family of curves $\{ \Phi_{j}(x) \}$ as a \emph{minimizing sequence} if $\lim_{j \rightarrow \infty} \sup_{x \in X} |\Phi_{j}(x)| = L_{i}$. If $|\Phi_{j}(x_{j})|$ converges to $L_{i}$ for some sequence of parameters $\{ x_{j}\}$ where $x_{j} \in X$ and $\{ \Phi_{j}(x) \}$ is a minimizing sequence, then let us call $\Phi_{j}(x_{j})$ as a \emph{min-max sequence}. Also let the \emph{critical set} $\Lambda(\{ \Phi_{j} \})$ be a set of stationary varifolds can be obtained by the limit of min-max sequence induced by  $\{ \Phi_{j}(x) \}$.

Let the set of critial geodesic $W_{L_{i}}$ be a set of stationary varifolds whose support is a simple closed geodesic and length is $L_{i}$. Moreover, We denote $W_{L_{i},j}$ and $W^{j}_{L_{i}}$ by the elements in $W_{L_{i}}$ whose support has Morse index less than or equal to $j$ or larger than or equal to $j$, respectively.

By curve shortening flow argument in \cite{Gr}, we have the following version of pull-tight type theorem (cf. Lemma 8.1 of \cite{Gr}).
\begin{thm}
Suppose $(S^{2},g)$ is endowed with a smooth metric. For any minimizing sequence $\{ \Phi_{j} \}$ of $i$-sweepouts, there is a deformed minimizing sequence $\{\hat{ \Phi}_{j} \}$ of $\{ \Phi_{j} \}$ satisfying the following property. For any $s>0$, there is some $0<a<L_{i}$ satisfying
    \begin{equation}
        \{ \hat{\Phi}_{j}(x) \in IV_{1}(S^{2}) : |\hat{\Phi}_{j}(x)| \ge L_{i}-a \} \subset \bigcup_{\gamma \in \Lambda(\{ \Phi_{j} \})\cap W_{L_{i}}} B^{F}_{s}(\gamma)
        \end{equation}
for all sufficiently large $j$, where $B^{F}_{s}(\gamma)$ is a $F$-metric ball with center $\gamma$. Moreover, the multiplicity of geodesics in the critical set is $1$.        
\begin{proof}
By the compactness theorem of simple closed geodesics (Theorem A.1), for $V \in V_{1}(S^{2})$, $F(V,W_{L_{i}})$ is well-defined. 

Given a simple closed curve $\alpha$, denote $H(t,\alpha)$ as a curve deformed by curve shortening flow at time $t$. Note that $H(t,\alpha)$ is an embedded curve for any $t \ge 0$. We set $\hat{\Phi}_{j}(x):= H(t_{j},\Phi_{j}(x))$ for $x \in X$ where $\{ t_{j} \}$ is a sequence such that $t_{j} \rightarrow \infty$ as $j \rightarrow \infty$.

We argue by contradiction and assume that there is a sequence of varifolds $\{ \hat{\Phi}_{j_{k}}(x_{j_{k}}) \}_{k \in \mathbb{N}}$ such that $|\hat{\Phi}_{j_{k}}(x_{j_{k}})| \ge L_{i}-a_{k}$ and $F(\hat{\Phi}_{j_{k}}(x_{j_{k}}),W_{L_{i}}) \ge s$ and $a_{k} \rightarrow 0$, $j_{k} \rightarrow \infty$ as $k \rightarrow \infty$. After passing to a subsequence, $\lim_{k \rightarrow \infty}|\hat{\Phi}_{j_{k}}(x_{j_{k}})| = L_{i}$ and $\lim_{k \rightarrow \infty} F(\hat{\Phi}_{j_{k}}(x_{j_{k}}),W_{L_{i}}) \ge s$. However, as time $t_{j_{k}}$ goes to infinity, $\hat{\Phi}_{j_{k}}(x_{j_{k}}) = H(t_{j_{k}},\Phi_{j_{k}}(x_{j_{k}}))$ converges to geodesics in $W_{L_{i}}$ or simple closed geodesics whose length is smaller than $L_{i}$. This contradicts to our assumption.  

Since curves converge to a multiplicity one geodesic or a single point by the curve shortening flow at limit time, the multiplicity of geodesics in the critical set is $1$.
\end{proof}
\end{thm}

By the proceeding pull-tight argument with curve shortening flow, we can construct simple closed geodesics for each $i$-sweepout. Let us consider spheres with bumpy metric. From the bumpiness of the metric $(S^{2},g)$ and Corollary 8.3 of \cite{Gr} (from Lusternik-Schnirelmann's topological argument), we obtain the existence of three simple closed geodesics with distinct lengths as following.

\begin{cor}
There are at least three simple closed geodesics on the bumpy sphere $(S^{2},g)$ with length $L_{1}, L_{2}$ and $L_{3}$.
\begin{proof}
Assume that two critical lengths from different sweepouts are the same. Corollary 8.3 of \cite{Gr} asserts that there are infinitely many simple closed geodesics if two critical length are the same. But it does not happen in the bumpy sphere $(S^{2},g)$ by Corollary A.2. Hence, the lengths of three simple closed geodesics are distinct each other. The pull-tight argument of Theorem 2.1 by flow gives three simple closed geodesics with length $L_{1},L_{2}$ and $L_{3}$.
\end{proof}
\end{cor}
Next we prove pull-tight type theorem to avoid the geodesics with large index in the sweepout. We need this theorem to prove the lower index bound in the proof of Theorem 1.2. This theorem is an analog of Theorem 4.9 of \cite{MN3} and is deduced by the constructive proof of Deformation theorem A of \cite{MN2}. We use the notation in Theorem 5.1 of \cite{MN2} in our proof.

Note that Deformation Theorem A in \cite{MN2} can be applied to sweepouts with smooth topology (See Remark 1.3 of \cite{MN2}). To be more precise, all deformations in the theorem are isotopies, and we can get smooth deformations by mollifying continuous deformations.

\begin{thm}
Suppose $(S^{2},g)$ is a $2$-sphere with a bumpy metric. For any minimizing sequence $\{ \Phi_{j} \}$ which is an $i$-sweepout, there is a deformed minimizing sequence $\{\hat{ \Phi}_{j} \}$ of $\{ \Phi_{j} \}$ satisfying the following property. For any small $s>0$, there is some $0<a<L_{i}$ satisfying
    \begin{equation}
        \{ \hat{\Phi}_{j}(x) \in IV_{1}(S^{2}) : |\hat{\Phi}_{j}(x)| \ge L_{i}-a \} \in \bigcup_{\gamma \in  \Lambda(\{ \Phi_{j} \}) \cap W_{L_{i},i}} B^{F}_{s}(\gamma)
        \end{equation}
for all sufficiently large $j$. Moreover, the multiplicity of geodesics in the critical set is $1$.         
\begin{proof}
By applying Theorem 2.1 and for minimizing sequence $\{ \Phi_{j} \}$, there is a deformed minimizing sequence $\{ \Psi_{j} \}$ such that for any $s>0$ there is some $0<a'<L_{i}$ satisfying (4):
\begin{equation*}
\{ \Psi_{j}(x) \in IV_{1}(S^{2}) : |\Psi_{j}(x)| \ge L_{i}-a' \} \in \bigcup_{\gamma \in W_{i}} B^{F}_{s}(\gamma)
\end{equation*}
for all sufficiently large $j$. 

Let $0<s<\epsilon$ where $\epsilon$ is a constant from Deformation Theorem A of \cite{MN2}. We deform $\{ \Psi_{j} \}$ to $\{ \hat{\Phi}_{j} \}$ by the deformation of Deformation Theorem A. Notice that $|\Psi_{j}(x)|<L_{i}-a'$ if $x \in X \setminus U_{i,2 \epsilon}$ by the above and $s< 2 \epsilon$. We divide into three cases and show that the curve with length close to $L_{i}$ is not in the neighborhood of simple closed geodesics with large index. 

If $x \in X \setminus U_{i,2 \epsilon}$, then $\hat{\Phi}_{i}(x) = \Psi_{i}(x)$ and $F(\gamma, \hat{\Phi}_{i}(x))>\epsilon>s$ where $Index(\gamma) \ge i+1$.

If $x \in U_{i,2 \epsilon} \setminus U_{i,5 \epsilon/4}$, we know that
\begin{equation}
    |\hat{\Phi}_{i}(x)| \le |(F_{H_{i}(x,1)})_{\sharp} (\Psi_{i}(x))|
\end{equation}
and
\begin{equation}
    \lim_{i \rightarrow \infty} \sup_{x \in X} ||\Psi_{i}(x)|-|(F_{H_{i}(x,1)})_{\sharp} (\Psi_{i}(x))||=0
\end{equation}
where $H_{i}(x,t)$ is a perturbation homotopy of Deformation Theorem A. Then (6) and (7) imply that
\begin{equation*}
    |\hat{\Phi}_{i}(x)| < L_{i} - \frac{a'}{2}
\end{equation*}
for large $i$.

We consider the last case of $x \in U_{i,5 \epsilon /4}$. By the deformation we obtain the following bound of the length:
\begin{equation*}
    |\hat{\Phi}_{i}(x)| \le L_{i} - \frac{c_{0}}{20}.
\end{equation*}
Now we take $a = \min (\frac{a'}{2},\frac{c_{0}}{20})$ and we obtain our claim by combining all cases. Multiplicity one property is not changed by our deformation so that the multiplicity of geodesics in critical set is also $1$.
\end{proof}
\end{thm}
\section{$F$-distance and Hausdorff distance of 1-dimensional integral varifolds}

In this section, we introduce a useful lemma to deal with $1$-dimensional integral varifolds induced by simple closed curves. We prove that if two integral varifolds induced by simple closed curves are sufficiently close in $F$-metric, then two varifolds are close in Hausdorff distance sense.

We show the following lemma by obtaining the bound of Hausdorff measure of the part of one varifold lying on the outside of small tubular neighborhood of another varifold. Since the support of each varifold is connected, every point on one varifold cannot go very far from another varifold by this estimate. This argument is only valid for $1$-dimensional varifolds. Let us denote Hausdorff distance between two nonempty subsets $A$ and $B$ of $(S^{2},g)$ as
\begin{equation*}
    d_{\mathcal{H}}(A,B) = \max \Big\{ \sup_{a \in A} d(a,B), \sup_{b \in B} d(b,A) \Big\},
\end{equation*}
where $d(a,B) = \inf _{b \in B} d(a,b)$ and $d(a,b)$ is an (intrinsic) distance between two points on $(S^{2},g)$.
\begin{lem}
Let $V, W \in IV_{1}(S^{2})$ be integral varifolds induced by simple closed curves. For $h>0$, if $F$-distance between two varifolds $V$ and $W$ satisfies
\begin{equation}
    F(V,W) < h^{2}/10,
\end{equation}
then 
\begin{equation*}
    d_{\mathcal{H}}(supp(V),supp(W)) < h.
\end{equation*}
\begin{proof}
Let us define the function $f: G_{1}(S^{2}) \rightarrow \mathbb{R}$ as
\begin{equation*}
f(x,\pi) = \max ( d(x,supp(W)),1),    
\end{equation*}
where $\pi \in T_{x}S^{2}$. Note that $f$ satisfies $f \ge 0$, $|f| \le 1$ and $Lip(f) \le 1$ so it satisfies the condition of test function for $F$-distance, and $f(x,\pi)=0$ at $x \in supp(W)$. Also, it suffices to consider the case of $|V \llcorner (S^{2} \setminus N_{\frac{h}{2}}(supp(W)))| \neq 0$. Now we obtain the estimate of the length of pieces of curves outside the tubular neighborhood $N_{\frac{h}{2}}(supp(W))$ of $supp(W)$ as

\begin{eqnarray}
    |V \llcorner (S^{2} \setminus N_{\frac{h}{2}}(supp(W)))| &\le& \Big( \frac{h}{2} \Big) ^{-1} \int f d(V \llcorner (S^{2} \setminus N_{\frac{h}{2}}(supp(W) ))\\ \nonumber
    &\le& \Big( \frac{h}{2} \Big) ^{-1} \int f dV \\ \nonumber
    &=& \Big( \frac{h}{2} \Big) ^{-1} \Big| \int f dV - \int f dW \Big| \\
    &\le& \Big( \frac{h}{2} \Big) ^{-1} F(V,W) < \Big( \frac{h}{2} \Big) ^{-1} \times \frac{h^{2}}{10}< \frac{h}{4}.
\end{eqnarray}
(9) is from $f \ge \frac{h}{2}$ on the outside of $N_{\frac{h}{2}}(supp(W))$, and (10) is from the $F$-distance assumption (8) and the definition of the $F$-distance (3).

Now we show that the whole support of $V$ is in an $h$-neighborhood of $supp(W)$ by using a triangle inequality. Note that for given $x \in supp(V) \cap (S^{2} \setminus N_{\frac{h}{2}}(supp(W) ))$ and $y \in supp(W)$, there exists $z_{x,y}$ satisfying $d(x,z_{x,y}) \le h/4$ and $d(y,z_{x,y}) = h/2$. Then we have
\begin{equation*}
d(x,y) \le d(x,z_{x,y}) +d(z_{x,y},y) \le h/4+ h/2 < h,
\end{equation*}
for any $x \in supp(V) \cap (S^{2} \setminus N_{\frac{h}{2}}(supp(W) ))$ and $y \in supp(W)$. Hence, the entire support of $V$ is in $N_{h}(supp(W))$. Equivalently, we obtain that the entire support of $W$ is in $N_{h}(supp(V))$ and have the upper bound of Hausdorff distance.
\end{proof}

\end{lem}
\section{Local Min-max Families}
In this section, we prove the $1$-varifold version of the local min-max theorem in \cite{MN3}, where the original version is Theorem 5 of \cite{W2}. While the proof is the simpler version of the proof of Theorem 6.1 of \cite{MN3}, we include the argument for the sake of completeness. Here we denote ${B}^{k}$ as an open unit ball in $\mathbb{R}^{k}$.
\begin{thm}
Let $\gamma$ be a simple closed geodesic with Morse index $k$ and multiplicity one. For every $\beta>0$, there is $\epsilon>0$ and a smooth family $\{ F_{v}\}_{v \in \bar{B}^{k}} \subset Diff(S^{2})$ such that 
\begin{enumerate}[label=(\roman*)]
\item $F_{0} = Id, \, F_{-v} = F_{v}^{-1}$ for all $v \in \bar{B}^{k}$; 
\item the function
\begin{equation*}
    L^{\gamma}: \bar{B}^{k} \rightarrow [0, \infty), \,\,\,\, L^{\gamma}(v) = |(F_{v})_{\sharp} \gamma|,
\end{equation*}
is strictly concave;
\item $||F_{v}-Id||_{C^{1}}< \beta$ for all $v \in \bar{B}^{k}$;
\end{enumerate}
and such that for every $V \in IV_{1}(S^{2})$ induced by a simple closed curve with $F(V, \gamma) < \epsilon$, we have
\begin{equation*}
\max_{v \in \bar{B}^{k}} |(F_{v})_{\sharp}V| \ge |\gamma|
\end{equation*}
with equality only if $\gamma = (F_{v})_{\sharp}V$ for some $v \in \bar{B}^{k}$.
\begin{proof}
In this proof, we adopt notations in \cite{MN3}. We define a smooth family $\{ F_{v} \}_{v \in \bar{B}^{k}} \subset Diff(S^{2})$ with properties satisfying (i) and (ii) constructed by the first $k$ linear combination of normal eigensections $\{ X_{i} = \frac{d}{dt} F_{te_{i}|t=0} \}_{1 \le i \le k}$ of the stability operator of $\gamma$ on $(S^{2},g)$ which are $L^{2}$-orthonormal each other. Let us define
\begin{equation*}
    P^{V} : \bar{B}^{k} \rightarrow \mathbb{R}^{k}, \,\,\,\, P^{V}(v) = \sum_{i=1}^{k} \Big( \int  \eta_{i} d ((F_{v})_{\sharp}V) \Big) e_{i},
\end{equation*}
where $e_{i}$'s are coordinate vectors in $\mathbb{R}^{k}$ and $\eta_{i}$ satisfies $\eta_{i} \le 1$, $\eta_{i}=0$ on $\gamma$ and $\nabla \eta_{i} = X_{i}$. By using the fact that $P^{V}$ is a regular diffeomorphism in the small neighborhood $B_{\delta}(0) \subset \bar{B}^{k}$ of the origin since $P^{V}(0)=0$ and $DP^{V}(0) = Id$, and our family satisfies (iii) by taking small $\delta>0$. Now let us consider the functional $L^{*}$ on $IV_{1}(S^{2})$ as
\begin{equation}
    L^{*}(V) = |V|+ (|\lambda_{1}|+1) \sum_{i=1}^{k} \Big(\int \eta_{i} dV \Big)^{2},
\end{equation}
where $\lambda_{1}$ is the first eigenvalue of the stability operator of $\gamma$. Then $\gamma$ is a strictly stable point of $L^{*}$. From Theorem 5 in \cite{W2} we obtain a tubular neighborhood $N_{h}(\gamma)$ such that $\gamma$ is a strict minimizer of $L^{*}$ in $N_{h}(\gamma)$. By setting $\epsilon=h^{2}/10$ and applying Lemma 3.1, we have $\gamma$ as a strict minimizer of the functional $L^{*}$ among varifolds induced by simple closed curves satisfying $F(V, \gamma)< \epsilon$.

We can choose $\delta '$ such that $F((F_{v})_{\sharp} \gamma,\gamma)< \epsilon/2$ for every $v \in \overline{B}^{k}_{\delta '}$ and argue by contradiction. Take $\{ F_{v} \}_{\overline{B}^{k}_{\delta '}}$ as a local min-max family and suppose there is a sequence of varifolds $\{ V_{j} \}$ induced by curves converging to $\gamma$ in the $F$-metric topology such that $(F_{v})_{\sharp} V_{j} \neq \gamma$ for all $v \in \overline{B}^{k}$ and
\begin{equation*}
\max_{v \in \overline{B}^{k}_{\delta '}} |(F_{v})_{\sharp}V_{j}| \le |\gamma|.    
\end{equation*}
Then $P^{V_{j}}|_{\overline{B}^{k}_{\delta'}}$ uniformly converges to $P^{\gamma}|_{\overline{B}^{k}_{\delta'}}$. Then there is $v_{j} \in \overline{B}^{k}_{\delta'}$ such that $P^{V_{j}}(v_{j}) =0$ by a degree argument. Then
\begin{equation*}
|(F_{v_{j}})_{\sharp} V_{j}| = L^{*}((F_{v_{j}})_{\sharp} V_{j}) >L^{*}(\gamma) = |\gamma|,    
\end{equation*}
for large $j$ and this gives contradiction. Hence, we obtained the desired conclusion.
\end{proof}
\end{thm}

\section{Squeezing a Family of Curves}
In this section, we prove a topological squeezing lemma to ensure the contractibility of a family of curves near a fixed simple closed geodesic. More precisely, this lemma asserts that the finite parameter family of simple closed curves within a small neighborhood of a simple closed geodesic can be squeezed into the simple closed geodesic. Hence, the family is nullhomotopic endowed with smooth topology and $F$-metric. The lemma is a varifold $F$-metric analog in smooth topology of Theorem 8.2 of \cite{Al} and Proposition 3.5 of \cite{MN1}, the homotopy lemma with flat topology. We apply curve shortening flow to squeeze a family of curves within a small neighborhood of a geodesic with smooth topology instead of adopting the construction of the homotopy from the work of Almgren \cite{Al}.

We construct the case of $\gamma$ is strictly stable with ambient negative Gaussian curvature first, and deal with general geodesics by metric perturbation. For the family near the strictly stable simple closed geodesic with ambient negative curvature, we flow the curves within the tubular neighborhood which can be foliated by mean convex curves to make all curves graphical. Then we apply a squeezing map of graphs to squeeze all curves in the family to the single geodesic. For the general case, we slightly perturb the metric near the geodesic to make it strictly stable with ambient negative curvature, and take a pullback homotopy in the smaller neighborhood.

Note that the squeezing lemma also yields the contractibility of finite parameter family of simple closed curves which are homologous to the center circle in any cylinder. Even though we construct the explicit homotopy via curve shortening flow in the small tubular neighborhood of geodesic, there is a diffeomorphism between any topological cylinder and the small tubular neighborhood so that we can extend this topological result to any cylinder.  This fact also follows from Smale's theorem on diffeomorphism group of $2$-sphere \cite{Sm}.

\subsection{Homotopy type of the curve near the geodesic in the tubular neighborhood} Suppose we have a simple closed curve sufficiently close to simple closed geodesic $\gamma$ in $F$-distance sense. We prove the proposition that the curve is homologous to $\gamma$ in the small tubular neighborhood $N_{h}(\gamma)$ of $\gamma$. We rule out the case of the nullhomotopic curve in the small tubular neighborhood to avoid the case of the curve converges to a single point by the curve shortening flow.

From now on, we adopt Fermi coordinates $c: [0,L] \times (-h,h) \rightarrow S^{2}$ on $N_{h}(\gamma)$ on the tubular neighborhood $N_{h}(\gamma)$ of the fixed geodesic $\gamma$ as in Appendix of \cite{G}. The metric on the tubular neighborhood of $\gamma$ is 
\begin{equation}
ds^{2} = J(x,y)^{2}dx^{2} + dy^{2}
\end{equation}
such that
\begin{enumerate} 
\item $J_{yy} = -KJ$
\item $J(x,0)=1$
\item $J_{x}(x,0) = \kappa =0$
\item $J_{y}(x,0) = 0$
\end{enumerate}
where $K$ is a Gaussian curvature and $\kappa$ is a geodesic curvature of $\gamma$.

\begin{prop}
If $\alpha$ is an integral varifold induced by a simple closed curve and $F(\alpha,\gamma)<h^{2}/10$, then $\alpha$ is homologous to $\gamma$ in a small tubular neighborhood $N_{h}(\gamma)$ of $\gamma$.
\begin{proof}
Let us assume that $\alpha$ is not homologous to $\gamma$ i.e. the curve is nullhomotopic in $N_{h}(\gamma)$ and show this lemma by contradiction. Note that the Hausdorff distance between two curves satisfies $d_{\mathcal{H}}(\alpha,\gamma)<h$ since $F(\alpha, \gamma)<h^{2}/10$ by Lemma 3.1.

Let the length of $\gamma$ be $L$. Throughout this proof we endow Fermi coordinate $c: [0,L] \times (-h,h) \rightarrow S^{2}$ on $N_{h}(\gamma)$, where the metric on the tubular neighborhood is given as (12). We lift this Fermi coordinate to the universal cover of the cylinder $\mathbb{R} \times (-h,h)$ and consider the lifted coordinate $C: \mathbb{R} \times (-h,h) \rightarrow S^{2}$, where $\pi \circ C = c$ and the lifted metric is given by $d \tilde{s}^{2} = J^{2} d \tilde{x}^{2} +dy^{2}$. Then we define the \textit{horizontal width} $W_{\sigma}$ of a closed curve $\sigma$ in $N_{h}$ as the maximum difference of $\tilde{x}$-coordinate on the lift:
\begin{equation*}
    W_{\sigma}:= \sup_{p_{1}, p_{2} \in \tilde{{\sigma}}} |\tilde{x}(p_{1})-\tilde{x}(p_{2})|,
\end{equation*}
where $C(\tilde{x}(z),y(z)):=z$ for every $z \in \sigma$. We point out that the lift $\tilde{\alpha}$ of $\alpha$ is a closed curve since $\alpha$ is nullhomotopic.  

We show that the horizontal width of $\alpha$ is less than $\frac{3L}{4}$. Suppose $p_{1}$ and $p_{2}$ are two points on $\tilde{\alpha}$ realizing the horizontal width of $\alpha$ and these two points divide $\tilde{\alpha}$ into two pieces. Let the lengths be $l_{1}$ and $l_{2}$. Let us define $l$ as a distance between $p_{1}$ and $p_{2}$. Then we have a following estimate:
\begin{equation*}
l \le \frac{l_{1}+l_{2}}{2} = \frac{|\alpha|}{2} < \frac{L+F(\alpha,\gamma)}{2}  = \frac{L+\delta}{2} < \frac{L+h^{2}}{2}. 
\end{equation*}
For the second inequality, we used $F(V,W) \ge ||V|-|W||$ for varifolds $V$ and $W$ by (3). Then by the triangle inequality, we have
\begin{eqnarray}
W_{\alpha} = |x(p_{1})-x(p_{2})| &\le& d(C(p_{1},0),p_{1})+d(p_{1},p_{2}) + d(p_{2},C(p_{2},0)) \nonumber \\
&\le& h+ \frac{L+h^{2}}{2} +h \nonumber \\
&<& \frac{3L}{4}
\end{eqnarray}
(13) comes from $h \ll L$. Since the horizontal width $W_{\alpha}$ of $\alpha$ satisfies $W_{\alpha}<\frac{3L}{4}$, there is a strip $\pi(C([a,a+L] \times (-h,h))) \subset N_{h}(\gamma)$ which does not contain any point of $\alpha$. Hence there exists a point $x \in \gamma$ such that $\alpha \cap B_{h}(x) = \emptyset$ so that $d_{\mathcal{H}}(\alpha,\gamma) \ge h$ and it contradicts $F(\alpha, \gamma)<h^{2}/10$ of Lemma 3.1. 
\end{proof}
\end{prop}
\subsection{Perturbation of the metric} In this subsection, we conformally deform the metric to make a geodesic to a strictly stable one with negative ambient Gaussian curvature. From this perturbation, we can construct the homotopy to deal with the geodesic by applying the pullback homotopy from the map will be constructed in the next subsection for the strictly stable and negative ambient curvature case.

 Now we perturb the metric near $N_{h}(\gamma)$ to convert $\gamma$ to a strictly stable geodesic with negative ambient Gaussian curvature. Let $M> \max K_{(S^{2},g)}$ and let $\nu$ be a unit normal vector field on $\gamma$, and for later estimates, let us take small $h$ satisfying
 \begin{equation}
    h<\max(M^{-1},1/10).
 \end{equation}
\begin{prop}
Given a simple closed geodesic $\gamma$, for any $\beta>0$, there exists a smooth bump function $\phi_{\beta}:S^{2} \rightarrow \mathbb{R}$ such that
\begin{enumerate}[label=(\roman*)]
\item $\phi_{\beta}(x)=0$ when $x \notin N_{h}(\gamma)$
\item $\phi_{\beta}(x)<0$ for $x \in \gamma$
\item $||\exp(2\phi_{\beta}) - 1||_{C^{0}} < \beta$
\item $\nabla_{\dot{\gamma}} \phi_{\beta} =0$ on $\gamma$ [$\phi_{\beta}$ is constant on $\gamma$]
\item The Hessian $\partial^{2} \phi_{\beta} (\nu,\nu)=M$ on $\gamma$.
\end{enumerate}
\begin{proof}
Define $\phi:S^{2} \rightarrow \mathbb{R}$ in terms of its Fermi coordinate $c: [0,L] \times (-h,h) \rightarrow S^{2}$ on $N_{h}(\gamma)$ in (12) as
\begin{equation}
    \phi(c(x,y)) = - A e^{-B/(h^{2}-y^{2})}
\end{equation}
for $x \in [0,L]$ and $y \in (-h,h)$, and put $\phi(z)=0$ if $z \notin N_{h}(\gamma)$. This guarantees (i), (ii) and (iv). Note that $\phi$ is smooth on $S^{2}$. We can calculate $\partial^{2} \phi (\nu,\nu)$ on $\gamma$ as
\begin{equation*}
\partial^{2} \phi(\nu,\nu)|_{z = \phi(c(x,0))} = \partial^{2} \phi(c(x,0)) /\partial y ^{2} =2AB h^{-4} e^{-B/h^{2}}  = -2 B h^{-4} \phi(c(x,0)).
\end{equation*}
For any $\epsilon>0$, let us take large $B$ by
\begin{equation}
B = h^4 \epsilon^{-1}M,
\end{equation}
and take $A = \epsilon e^{B/h^{2}}/2$ so that $\phi(c(x,0)) = -\epsilon/2$. Then we set that $\partial^{2} \phi(c(x,0)) (\nu, \nu)=M$ and $||\phi||_{C^{0}}< \epsilon$. Thus for any $\beta>0$ we can construct a bump function $\phi_{\beta}$ to satisfy $||\exp(2\phi_{\beta}) - 1||_{C^{0}} < \beta$ by taking appropriate $\epsilon$ and corresponding $A$ and $B$ in (15). This verifies (iii) and (v) by our choice of $A$ and $B$. 
\end{proof}
\end{prop}
For a fixed $\beta>0$, we have $\phi_{\beta}$ by Proposition 5.2 and the metric $g_{\beta}=\exp(2 \phi_{\beta}) g$ in Fermi coordinate is
\begin{equation}
ds_{\beta}^{2} = J(x,y)^{2}\exp(2 \phi_{\beta}) dx^{2} +\exp (2 \phi_{\beta}) dy^{2}
\end{equation}
and we have a canonical diffeomorphism $G_{\beta}: (S^{2},g) \rightarrow (S^{2},g_{\beta})$ given by $G_{\beta}(z) \equiv z$ for all $z \in S^{2}$. Let $\gamma_{g_{\beta}} := G_{\beta}(\gamma)$ and $N_{h,g_{\beta}}(\gamma_{g_{\beta}}) : = G_{\beta}(N_{h}(\gamma_{g_{\beta}}))$ in $(S^{2},g_{\beta})$. 
\begin{lem} 
For any $\beta>0$, $\gamma_{g_{\beta}}$ is a strictly stable geodesic in $N_{h,g_{\beta}}(\gamma_{g_{\beta}})$. Moreover, the Gaussian curvature satisfies $K_{g_{\beta}}(x)<0$ for $x \in \gamma_{g_{\beta}}$ on $(S^{2},g_{\beta})$.
\begin{proof}
By the calculation of the change of second fundamental form by conformal deformation in Besse \cite{Be}, we obtain the geodesic curvature as
\begin{equation}
    \kappa_{g_{\beta}} = e^{-\phi}\Big(\kappa_{g}-\frac{\partial \phi_{\beta}}{\partial \nu}\Big),
\end{equation}
where $\partial \phi/\partial \nu$ is a normal derivative of $\phi$ to $\gamma$ with respect to the metric $g_{\beta}$ so that the geodesic curvature on $\gamma_{g_{\beta}}$ becomes $0$ and so $\gamma_{g_{\beta}}$ is a geodesic.

We prove that the Gaussian curvature of $(S^{2},g_{\beta})$ on $\gamma_{g_{\beta}}$ is negative and $\gamma_{g_{\beta}}$ is a stable geodesic in $N_{h,g_{\beta}}(\gamma_{g_{\beta}})$.

The Gaussian curvature $K_{g_{\beta}}(x)$ at $x \in \gamma_{g_{\beta}}$ is
\begin{eqnarray}
\nonumber K_{g_{\beta}}(x)&=& \nonumber \partial_{2} \Gamma^{2}_{11} - \partial_{1} \Gamma^{2}_{21} + \Gamma^{1}_{11} \Gamma^{2}_{21} + \Gamma^{1}_{12} \Gamma^{2}_{22} - \Gamma^{1}_{21} \Gamma^{2}_{11}- \Gamma^{2}_{21} \Gamma^{2}_{12} \\
&=& K(x) - \partial^{2} \phi_{\beta} (\nu,\nu) = K(x)-M<0.
\end{eqnarray}
The inequality in (19) comes from $M >\max(\max K_{(S^{2},g)},0 )$. Suppose $f:\gamma \times (-\epsilon,\epsilon) \rightarrow \mathbb{R}$ to be a normal variation on $\gamma_{g_{\beta}}$. Then since the stability operator of $f$ is $\delta_{\gamma_{g_{\beta}}}(f) = \int_{\gamma_{g_{\beta}}}(|\nabla_{\gamma_{g_{\beta}}} f|^{2} -K_{g_{\beta}}f^{2})ds$ and is positive by $K_{g_{\beta}}<0$. We conclude that $\gamma_{g_{\beta}}$ is a strictly stable geodesic.
\end{proof}
\end{lem}

\subsection{Squeezing lemma} Now we prove the following squeezing lemma by constructing the homotopy between family of curves. Our construction mainly relies on the curve shortening flow and the squeezing map in the sense of Fermi coordinate. Our explicit construction of the squeezing map will be used to settle the quantitative interpolation lemmas in Section 6. 

For a strictly stable simple closed geodesic $\gamma$ with negative ambient Gaussian curvature, some tubular neighborhood $N_{h}(\gamma)$ of $\gamma$ can be foliated by simple closed curves $\{ \partial N_{t}(\gamma) \}_{t \in [0,h]}$ whose curvature vector points toward $\gamma$. To see this, we adopt the Fermi coordinate $c:[0,L] \times (-h,h) \rightarrow S^{2}$ on $N_{h}(\gamma)$ in (12) and follow the calculation in Appendix of \cite{G}. Denote the orthonormal frame field in $N_{h}(\gamma)$ by $e_{1} = (1/J) c_{*}(\partial/\partial x)$ and $e_{2} = c_{*}(\partial /\partial y)$, and define 
\begin{equation}
\gamma_{t}:= c( \{y=t \})
\end{equation}
as a level curve for $t \in [-h,h]$. Then (A.1) in \cite{G} gives that the geodesic curvature $\kappa_{t}(x)$ (the positive direction is $e_{2}$) of the curve $\gamma_{t}$ at $c(x,t) \in \gamma_{t} \subset S^{2}$ to be:
\begin{equation*}
    \kappa_{t}(x) = -J_{y}(x,t)/J(x,t).
\end{equation*}
From the conditions (1),(2),(4) for (12), we have $J_{y}(x,t)/t>0$ for $0<|t|<h'$ for any $x \in [0,L]$ and small $h'>0$. This gives the mean convex foliation $\{ \partial N_{t}(\gamma) \}_{t \in [0,h']}$ of $\gamma$ (See Proposition 5.7 of \cite{MN2} for higher dimensional case with the first eigenfunction of Jacobi operator). Let us replace $h'$ by $h$ and call $N_{h}(\gamma)$ as a \emph{mean convex neighborhood}.

\begin{lem}
Let $(S^{2},g)$ be a sphere with bumpy metric, $\gamma$ be a simple closed geodesic, and $X$ be a simplicial complex with finite dimension $k$. There exists $\delta_{0} = \delta_{0}((S^{2},g), \gamma)>0$ with the following property:

For $0<\delta<\delta_{0}$, if $\Phi : X \rightarrow \mathcal{S}$ is a continuous map in the smooth topology such that
\begin{equation*}
    \sup \{ F(\Phi(x), \gamma) : x \in X \} < \delta,
\end{equation*}
then there is a homotopy $H:[0,1] \times X \rightarrow \mathcal{S}$ such that $H(0,x) = \Phi(x)$ and $H(1,x) = \gamma$ so that $\Phi$ is nullhomotopic.
\begin{rmk}
The contractiblity of a finite dimensional family of simple closed curves in a topological cylinder also follows from Smale's theorem on diffeomorphism group of Riemannian $2$-spheres \cite{Sm}.
\end{rmk}
\begin{rmk}
We can construct the squeezing map in the strictly stable case without the ambient curvature condition directly rather than taking a pullback homotopy since there exists a mean convex foliation on the tubular neigborhood of the geodesic. Note that the curves in this foliation are not necessarily to be the level curves in (20). We need to distinguish these two possible foliations to obtain the quantitative $F$-distance in Section 6 with the technical reason. Hence, we divide the cases to the strictly stable geodesic with negative ambient curvature cases and the general geodesic cases in the proof of Lemma 5.4.
\end{rmk}
\begin{proof}
First we prove the lemma when $\gamma$ is a strictly stable geodesic and whose ambient negative Gaussian curvature is negative. We are able to suppose the only simple closed geodesic is in the tubular neighborhood $N_{h}(\gamma)$ of $\gamma$ since the metric is bumpy.

For $h$ such that $N_{h}(\gamma)$ is a mean convex neighborhood, we set $\delta_{0} = h^{2}/10$. Then the entire support of $\Phi(x)$ are in $N_{h}(\gamma)$ by Lemma 3.1 for any $x \in X$. We apply the curve shortening flow to the curve $\Phi(x)$ for each $x$. Note that the curve shortening flow deforms the family of curves continuously and each curve converges to a single point or a simple closed geodesic as time goes to infinity (See Proposition 1.4 and Theorem 3.1 of \cite{G}).

Since $\Phi(x)$ is homologous with $\gamma$ in $N_{h}(\gamma)$ by Proposition 5.1 and $N_{h}(\gamma)$ is a mean convex neighborhood, each $\Phi(x)$ uniformly $C^{\infty}$-converges to $\gamma$ as time goes to infinity by the avoidance principle of the curve shortening flow.

Let $Y(\cdot,t)|_{\Phi(x)}$ be the curve shortening flow for $\Phi(x)$ where $Y(\cdot,0)|_{\Phi(x)} = supp(\Phi(x))$ and $Y_{t}(y,t)|_{\Phi(x)} = \kappa N$ where $\kappa$ is a geodesic curvature and $N$ is a normal vector field on $\Phi(x)$. Then for any small $s>0$ there is a large time $t_{0}$ such that for $t>t_{0}$ and any $x \in X$, the total curvature $\int_{Y(\cdot,t)} |\kappa| < s$ since curves uniformly converge to $\gamma$. We construct the homotopy $H_{1} : [0,1] \times X \rightarrow IV_{1}(S^{2})$ induced by the flow as
\begin{equation*}
    H_{1}(t,x) = Y(\cdot,2t_{0}t)|_{\Phi(x)},
\end{equation*}
for all $x \in X$ and $0 \le t \le 1$. We can take $t_{0}$ to make each curve $H_{1}(t,x)$ to be graphical in $N_{h}(\gamma)$ by taking sufficiently small $s>0$ and applying Lemma A.2 of \cite{G}. 

We construct a second homotopy between curves $H_{1}(1,x)$ and $\gamma$ by using the squeezing map (For similar constructions, see \cite{KLS} and \cite{MN2}). Let $c:[0,L] \times (-h,h) \rightarrow (S^{2},g)$ be the Fermi coordinate system. Since each curve $H_{1}(1,x)$ is graphical, we can parametrize each curve $H_{1}(1,x)$ as $H_{1}(1,x) = \{c(y,g(x,y)) | y \in [0,L] \}$, where $g : X \times [0,L] \rightarrow [-h,h]$ is a parametrization of the height function part in Fermi coordinate of (12). We construct the second homotopy by
\begin{equation*}
    H_{2}(t,x) = \{ c(y,(1-t)g(x,y)) | y \in [0,L]\}
\end{equation*}
for $0 \le t \le 1$ and $H_{2}(1,x)= \gamma$ for every $x \in X$. By combining two homotopies $H_{1}$ and $H_{2}$ by
\begin{equation}
H=H_{1} \cdot H_{2},
\end{equation}
we prove that $\Phi$ is nullhomotopic if $\gamma$ is a strictly stable geodesic whose ambient Gaussian curvature is negative. 

Now we consider a general geodesic $\gamma$. Let us fix some $\beta>0$ and consider the perturbed metric $(S^{2},g_{\beta})$ as in (17) and corresponding diffeomorphism $G_{\beta}$. Then $\gamma_{g_{\beta}}$ is a strictly stable geodesic and there exists a tubular neighborhood $N_{h'}(\gamma_{g_{\beta}})$ which is a mean convex neighborhood of $\gamma_{g_{\beta}}$. Note that the tubular neighborhood $N_{h'}(\gamma_{g_{\beta}})$ here should be taken as smaller one than $N_{h,g_{\beta}}(\gamma_{g_{\beta}})$. Let us take $\delta_{0} = h'^{2}/20$. Then we have
\begin{align}
   \nonumber F(&(G_{\beta})_{\sharp}\Phi(x),(G_{\beta})_{\sharp}\gamma) \\ \nonumber &= \sup \Big\{ \Big| \int f d((G_{\beta})_{\sharp}\Phi(x)) - \int f d((G_{\beta})_{\sharp}\gamma) \Big| : f \in C_{c}(G_{1}(S^{2})), |f| \le 1 , Lip(f) \le 1 \Big\} \\
    &= \sup \Big\{ \Big| \int f JG_{\beta} d\Phi(x) - \int f JG_{\beta} d\gamma \Big| : f \in C_{c}(G_{1}(S^{2})), |f| \le 1 , Lip(f) \le 1 \Big\},
\end{align}
where $JG_{\beta}$ is a Jacobian of $G_{\beta}$. Here we estimate the bound of $|JG_{\beta}|$ and $Lip(JG_{\beta})$ on $G_{1}(N_{h}(\gamma))$:
\begin{align}
|f(JG_{\beta})| &\le |f| \le 1, \\ \nonumber Lip(f(JG_{\beta})) & \le \Big( Lip(f) \sup_{(x, \pi) \in G_{1}(N_{h}(\gamma))} |JG_{\beta}| +  Lip(JG_{\beta}) \sup_{(x, \pi) \in G_{1}(N_{h}(\gamma))} |f| \Big) \\ & \le \Big( 1+ \sup_{N_{h}(\gamma)} \Big| \frac{\partial \phi_{\beta}}{\partial \nu} \exp(\phi_{\beta}) \Big| \Big) \\  & \le \Bigg( 1  +  \sup_{y \in (-h,h)} \Bigg| \frac{2ABy e^{-\frac{B}{h^{2}-y^{2}}}}{(h^{2}-y^{2})^{2}} \exp(\phi_{\beta}) \Bigg| \Bigg) \\ &\le 1+Mh<2. 
\end{align}
(24) comes from $JG_{\beta} = \exp(\phi_{\beta}) \le 1$ and (25) comes from (15). We use the choice of $h$ in (14), the choice of $B$ in (16) and the decreasing property of $e^{-B/(h^{2}-y^{2})}/(h^{2}-y^{2})^{2}$ in $[0,h]$ for small $h$ and large $B$ for (26). By (22), (23), (26) and our choice of $\delta_{0}$, we obtain
\begin{align}
     \nonumber F(&(G_{\beta})_{\sharp}(\Phi(x)),(G_{\beta})_{\sharp}\gamma)\\ \nonumber
    &= \sup \Big\{ \Big| \int f JG_{\beta} d\Phi(x) - \int f JG_{\beta} d\gamma \Big| : f \in C_{c}(G_{1}(S^{2})), |f| \le 1 , Lip(f) \le 1 \Big\} \\ \nonumber &= 2\sup \Big\{ \Big| \int \frac{f JG_{\beta}}{2} d\Phi(x) - \int \frac{f JG_{\beta}}{2} d\gamma \Big| : f \in C_{c}(G_{1}(S^{2})), |f| \le 1 , Lip(f) \le 1 \Big\} \\ \nonumber
    &\le 2 F(\Phi(x),\gamma) \\ &< h'^{2}/10,
\end{align}
and this implies $(G_{\beta})_{\sharp} (\Phi(x)) \in N_{h'}(\gamma_{g_{\beta}})$ for any $x \in X$ by Lemma 3.1. We now obtain the homotopy $H'_{\beta}:[0,1] \times X \rightarrow IV_{1}(S^{2})$ given by (21), where the endowed metric of $S^{2}$ is $g_{\beta} = \exp(\phi_{\beta})g$. Then we take the pullback homotopy $H_{\beta}:[0,1] \times X \rightarrow IV_{1}(S^{2})$ in $(S^{2},g)$ by
\begin{equation}
H_{\beta} := G_{\beta}^{*}H'_{\beta}.    
\end{equation}
We obtained the desired homotopy $H_{\beta}$ in the general geodesic case.  
\end{proof}
\end{lem}
\section{Quantitative estimate of $F$-distance and the interpolation lemma}
In this section, we prove the interpolation lemma (Lemma 6.4) between two families of simple closed curves which are close each other based on the homotopy we constructed in Section 5. We construct an explicit interpolation homotopy whose lengths do not exceed the critical length along the homotopy by composing local min-max diffeomorphism from Theorem 4.1 and squeezing map in Lemma 5.4.

We proved the contractibility of a finite-parameter family of simple closed curves which are homotopic to the center circle in any cylinder in Lemma 5.4. Hence we know that two family of simple closed curves which are homologous to the geodesic in the small tubular neighborhood of the geodesic are homotopic. However, it does not directly give the information on the length or $F$-distance bound along the homotopy even in the small tubular neighborhood.
\subsection{F-distance estimate}

In this subsection, we prove the quantitative estimate of $F$-distance along the squeezing map inspired by Lemma 5.4. First we show the quantitative estimate of $F$-distance between a fixed strictly stable geodesic and a homologous simple closed curve within its small tubular neighborhood when the length of the curve is not too larger than that of the geodesic when the ambient curvature is negative.

We prove the $F$-distance estimate of simple closed curves whose proof is inspired by the proof of Quantitative Constancy theorem for stationary varifolds of Section 2 of \cite{SZ}. They proved that if the total mass of the codimension $1$ stationary varifold in a closed $(n+1)$-dimensional manifold is mostly concentrated in a tubular neighborhood of a $2$-sided, closed and embedded hypersurface, then the varifold distance estimate holds between the normalized varifold and the normalized hypersurface. In our geodesic setting, we show that if the simple closed curve is in a tubular neighborhood of strictly stable geodesic with ambient negative curvature and whose length is bounded above, then $F$-distance between the simple closed curve and the strictly stable geodesic satisfies the quantitative estimate.

Let $\gamma$ be a strictly stable simple closed geodesic in $(S^{2},g)$ and Gaussian curvature on $\gamma$ is negative. For Lemma 6.1 and Lemma 6.2, we take sufficiently small $h$ such that $N_{h}(\gamma)$ is a mean convex neighborhood. Let us take the Fermi coordinate $c: [0,L] \times (-h,h) \rightarrow S^{2}$ of (12). Also define the distance function $dist_{g}$ between two tangent lines with a fixed point in the Grassmannian manifold $G_{1}(S^{2})$ endowed with the induced metric by $(S^{2},g)$ on the Grassmannian manifold $G_{1}(S^{2})$. 
\begin{lem}
Let $(S^{2},g)$, $\gamma$, $c$ be as above and $K(z)<0$ for $z \in N_{h}(\gamma)$. There exists $C=C(L)>0$ satisfying the following property: For $0<\epsilon<h^{2}$, if a simple closed curve $\alpha$ homologous to $\gamma$ in $N_{h}(\gamma)$ satisfies $|\alpha|<L+\epsilon$, then
\begin{equation}
F(\alpha,\gamma)< C(L)(h+ \sqrt{\epsilon}).    
\end{equation}
\begin{proof}
We work in the Fermi coordinate. Parametrize $\alpha$ by arclength as $\alpha(s)=c(x(s),y(s))$. Thus $J^{2}(x(s),y(s))x'^{2}(s) + y'^{2}(s) =1$ for $s \in [0,|\alpha|]$. Indeed, we have
\begin{equation}
    \int_{0}^{|\alpha|} \sqrt{J^{2}(x(s),y(s))x'^{2}(s) + y'^{2}(s)} ds = |\alpha| < L+ \epsilon.
\end{equation}
By condition (1) and (2) of $J$ in (12), we have
\begin{equation}
    J(x,y) \ge 1
\end{equation}
on $N_{h}(\gamma)$. Since $\alpha$ is a simple closed curve homologous to $\gamma$, by (31) we have
\begin{equation}
L \le \int_{0}^{|\alpha|} |x'(s)| ds \le \int_{0}^{|\alpha|} J(x(s),y(s))|x'(s)| ds.
\end{equation}
By subtracting (32) from (30), we obtain
\begin{equation}
    \int_{0}^{|\alpha|} \sqrt{J^{2}(x(s),y(s))x'^{2}(s) + y'^{2}(s)} - J(x(s),y(s))|x'(s)| ds < \epsilon.
\end{equation}
Consider the foliations by level curves $\{ \gamma_{t} \}_{t \in [-h,h]}$ in (20), and denote $\gamma_{z}$ as a leaf in the foliation $\{\gamma_{t} \}_{t \in [-h,h]}$ containing $z \in N_{h}(\gamma)$. Then we define $\theta(s) \in [-\pi,\pi]$ to be a continuous angle function between the tangent vector $dc_{(x(s),y(s))}(\partial/\partial x) \in T_{\alpha(s)}\gamma_{\alpha(s)}$ and the tangent vector $\alpha'(s)\in T_{\alpha(s)}\alpha$. Then $\tan \theta(s) = y'(s)/Jx'(s)$. By combining (33) with $J^{2}(x(s),y(s))x'^{2}(s) + y'^{2}(s) =1$, we deduce
\begin{equation}
    \int_{\alpha} 1- \cos(\theta(s)) ds < \epsilon.
\end{equation}
We obtain an estimate on the total angle of $\alpha$:
\begin{equation}
    \int_{\alpha} |\theta(s)| ds \le (L+\epsilon)^{\frac{1}{2}}\Big( \int_{\alpha} |\theta(s)|^{2} ds \Big)^{\frac{1}{2}} \le 2 L^{\frac{1}{2}} \Big( \int_{\alpha} \pi^{2}(1- \cos\theta(s)) ds \Big)^{\frac{1}{2}} < 2 \pi L^{\frac{1}{2}} \sqrt{\epsilon},
\end{equation}
where the first inequality comes from the Hölder's inequality and the second inequality is from the inequality $\theta^{2} \le \pi^{2}(1-\cos \theta)$ for $\theta \in [-\pi,\pi]$, and the last inequality follows from (34).

Now we deduce the relation between the total angle of the curve $\alpha$ and the total of the distance function $dist_{g}$ we defined before the statement of Lemma 6.1. There exists a constant $C>0$ such that $dist_{g}(T_{z}\alpha,T_{z}\gamma_{z}) \le C |\theta(s)|$ at $z = \alpha(s)$. Thus We have
\begin{equation}
   \int_{\alpha} dist_{g}(T_{z}\alpha,T_{z}\gamma_{z})d\alpha(z) \le C \int_{\alpha} |\theta(s)| ds.
\end{equation}
By combining (35) and (36), there exists a constant $C>0$ satisfying:
\begin{equation}
    \int_{\alpha} dist_{g}(T_{z}\alpha,T_{z}\gamma_{z})d\alpha(z) \le C L^{\frac{1}{2}} \sqrt{\epsilon}
\end{equation}
Let us denote $\pi : N_{h}(\gamma) \rightarrow \gamma$ as the orthogonal projection map onto $\gamma$. Now we estimate $F$-distance between $\alpha$ and $\gamma$. By combining the definition of $F$-distance (3) and (37) we may estimate:
\begin{align}
\nonumber \bigg| \int_{\alpha} &f(z,T_{z}\alpha) d\alpha(z) - \int_{\gamma} f(p,T_{p}\gamma) d\gamma(p) \bigg| \\ \nonumber &\le \bigg| \int_{\alpha} (f(z,T_{z}\alpha) - f(z,T_{z}\gamma_{z})) d \alpha(z)\bigg| + \bigg| \int_{\alpha}f(z,T_{z}\gamma_{z}) d\alpha(z) - \int_{\gamma} f(p,T_{p}\gamma) d\gamma(p) \bigg| \\ \nonumber &\le \int_{\alpha} dist_{g}(T_{z}\alpha,T_{z}\gamma_{z})d\alpha(z) + \bigg| \int_{\alpha}f(z,T_{z}\gamma_{z}) d\alpha(z) - \int_{\gamma} f(p,T_{p}\gamma) d\gamma(p) \bigg| \\  &\le CL^{\frac{1}{2}}\sqrt{\epsilon} + \bigg| \int_{\alpha}f(z,T_{z}\gamma_{z}) d\alpha(z) - \int_{\gamma} f(p,T_{p}\gamma) d\gamma(p) \bigg|.
\end{align}
By applying (3), (37), (38) repeatedly, we obtain the estimate:
\begin{align}
   \nonumber \bigg| \int_{\alpha} f(z,T_{z}\alpha)& d\alpha(z) - \int_{\gamma} f(p,T_{p}\gamma) d\gamma(p) \bigg| \\ \nonumber &\le CL^{\frac{1}{2}}\sqrt{\epsilon} + \bigg| \int_{\alpha}f(z,T_{z}\gamma_{z}) d\alpha(z) - \int_{\gamma} f(p,T_{p}\gamma) d\gamma(p) \bigg| \\ \nonumber &
    \le CL^{\frac{1}{2}} \sqrt{\epsilon}+ CLh +  \bigg| \int_{\alpha}f(\pi(z),T_{\pi(z)}\gamma) d\alpha(z) - \int_{\gamma} f(p,T_{p}\gamma) d\gamma(p) \bigg| \\ & 
    \le CL^{\frac{1}{2}} \sqrt{\epsilon}+ CLh + \bigg| \int_{\alpha}f(\pi(z),T_{\pi(z)}\gamma) d\alpha(z) - \int_{\pi(\alpha)} f(\pi(z),T_{\pi(z)}\gamma) J(z)d\gamma(\pi(z)) \bigg| \nonumber \\
    & + \bigg| \int_{\pi(\alpha)} f(\pi(z),T_{\pi(z)}\gamma) J(z)d\gamma(\pi(z)) - \int_{\gamma} f(p,T_{p}\gamma) d\gamma(p) \bigg| \nonumber \\ &
    = CL^{\frac{1}{2}} \sqrt{\epsilon}+ CLh + \bigg| \int_{\alpha}f(\pi(z),T_{\pi(z)}\gamma)(1- \cos\theta(z) ) d\alpha(z)  \bigg| \nonumber \\
    & + \bigg| \int_{\pi(\alpha)} f(\pi(z),T_{\pi(z)}\gamma) J(z)d\gamma(\pi(z)) - \int_{\gamma} f(p,T_{p}\gamma) d\gamma(p) \bigg| \nonumber\\ &
    \le CL^{\frac{1}{2}} \sqrt{\epsilon}+ CLh + C L^{\frac{1}{2}} \sqrt{\epsilon} + C \epsilon \\ &
    \le C(L)(h+ \sqrt{\epsilon}). \nonumber
\end{align}
Note that for (39), we used (3), (34), $J \ge 1$, $|\pi(\alpha)|<L+\epsilon$, and the fact that $\pi(\alpha)$ spans $\gamma$ at least once, where $\alpha$ is homologous to $\gamma$, and estimated similarly with (38). 
\end{proof}
\end{lem}
Now we have the $F$-distance estimate along the squeezing map in Lemma 5.4 when the geodesic $\gamma$ is strictly stable and the ambient Gaussian curvature is negative:

\begin{lem} Let $\gamma$ be a strictly stable geodesic on $(S^{2},g)$, Gaussian curvature $K(z)<0$ for $z \in N_{h}(\gamma)$, and $X$ be a $k$-dimensional simplicial complex. Then there exists $C = C(|\gamma|)>0$ satisfying the following property: For $0<\epsilon<h^{2}/10$, if $\Phi:X \rightarrow \mathcal{S}$ is a continuous map in the smooth topology such that
\begin{equation*}
    \sup \{ F(\Phi(x), \gamma) : x \in X \} < \epsilon,
\end{equation*}
then there is a homotopy $H:[0,1] \times X \rightarrow \mathcal{S}$ such that $H(0,x) = \Phi(x)$, $H(1,x) = \gamma$ and the following $F$-distance estimate holds along $H$:
\begin{equation*}
    \sup \{ F(H(t,x),\gamma) : x \in X \text{ and } t \in [0,1] \} < C(|\gamma|) \sqrt{\epsilon}.
\end{equation*}
\begin{proof}
Note that $N_{h}(\gamma)$ is a mean convex neighborhood and $K(z)<0$ for $z \in N_{h}(\gamma)$. Then $supp(\Phi(x)) \subset N_{\sqrt{10\epsilon}}(\gamma) \subset N_{h}(\gamma)$ for any $x \in X$ by Lemma 3.1 and $|\Phi(x)|< L+ \epsilon$ by (3). Let us consider the homotopy $H=H_{1} \cdot H_{2}$ in (21).

For the squeezing homotopy $H_{2}(t,x)= H(1/2+ t/2,x)$,
\begin{align}
    \int_{H(x,t)} ds &= \int_{0}^{|H(x,\frac{1}{2})|} \sqrt{J^{2}(x(s),(-2t+2)^{2}y(s))x'^{2}(s) + (-2t+2)^{2} y'^{2}(s)} ds \\
    &\le \int_{0}^{|H(x,\frac{1}{2})|} \sqrt{J^{2}(x(s),y(s))x'^{2}(s) + y'^{2}(s)} ds \\
     &= \nonumber \int_{H(x,\frac{1}{2})} ds,
\end{align}
for $1/2 \le t \le 1$ and where the parametrization in (40) of $H(t,x)$ is by the arclength of $H(\frac{1}{2},x)$. (41) comes from the convexity of $J$ and the fact that $J$ has a unique minimum on $y=0$ for a fixed $x$-coordinate.

Since the length functional on the homotopy by the curve shortening flow and squeezing map is monotonically decreasing, $|H(t,x)|<L+\epsilon$ for any $x \in X$ and $t \in [0,1]$. Moreover, $supp(H(t,x)) \subset N_{\sqrt{10\epsilon}}(\gamma)$ by the avoidance principle and the construction of mean convex foliation. By applying Lemma 6.1, we obtain $F(H(t,x),\gamma)<C(|\gamma|)\sqrt{\epsilon}$ for any $x \in X$ and $t \in [0,1]$. 
\end{proof}
\end{lem}
We now construct a squeezing homotopy with $F$-distance bound based on Lemma 5.4 in a sufficiently small neighborhood for general geodesics.

\begin{lem} Let $\gamma$ be a geodesic on $(S^{2},g)$, and $X$ be a $k$-dimensional simplicial complex. There exists $C = C(|\gamma|)>0$ and $\epsilon_{0}= \epsilon_{0}(\gamma)>0$ such that satisfying the following property: For $0<\epsilon<\epsilon_{0}$, if $\Phi : X \rightarrow \mathcal{S}$ is a continuous map in the smooth topology such that
\begin{equation*}
    \sup \{ F(\Phi(x), \gamma) : x \in X \} < \epsilon,
\end{equation*}
then there is a homotopy $H:[0,1] \times X \rightarrow \mathcal{S}$ such that $H(0,x) = \Phi(x)$, $H(1,x) = \gamma$ and the following $F$-distance estimate holds along $H$:
\begin{equation}
    \sup \{ F(H(t,x),\gamma) : x \in X \text{ and } t \in [0,1] \} < C(|\gamma|) \sqrt{\epsilon}.
\end{equation}
\begin{proof}
For a given $\beta$, we can take the conformal deformation $G_{\beta}$ of the metric $g_{\beta} = \exp(2 \phi_{\beta}) g$, whose metric is given by (17). Note that $\phi$ is in the form of $\phi(c(x,y)) = - A e^{-B/(h^{2}-y^{2})}$ by (15). We will specify $\beta$ later in this proof. Then we consider the homotopy $H$ in (28) for this given $\beta>0$ for the general geodesic $\gamma$ case. 

In a similar way to (22), we have
\begin{align*}
    F(&H(t,x),\gamma) = F((G_{\beta})^{\sharp}H'(t,x),(G_{\beta})^{\sharp}((G_{\beta})_{\sharp}\gamma))\\ &= \sup \Big\{ \Big| \int f d((G_{\beta})^{\sharp}H'(t,x)) - \int f d((G_{\beta})^{\sharp}((G_{\beta})_{\sharp}\gamma)) \Big| : f \in C_{c}(G_{1}(S^{2})), |f| \le 1 , Lip(f) \le 1 \Big\} \\
    &= \sup \Big\{ \Big| \int f (JG_{\beta})^{-1} dH'(t,x) - \int f (JG_{\beta})^{-1} d(G_{\beta})_{\sharp}\gamma \Big| : f \in C_{c}(G_{1}(S^{2})), |f| \le 1 , Lip(f) \le 1 \Big\}.
\end{align*}
Note that $(G_{\beta})_{\sharp}\gamma$ is a strictly stable geodesic, the ambient Gaussian curvature is negative on $(G_{\beta})_{\sharp}\gamma$.

Note that since $|JG_{\beta}| = \exp(\phi_{\beta}) \ge 1-\beta$,
\begin{equation}
    |f(JG_{\beta})^{-1}| \le (1-\beta)^{-1}|f|
\end{equation}
on $G_{1}(N_{h,g_{\beta}}(\gamma_{g_{\beta}}))$.
Moreover, by calculating similarly to (24)-(26), we have
\begin{align}
    \nonumber Lip(f(JG_{\beta})^{-1}) & \le \exp(-\phi_{\beta}) \Big( Lip(f) \sup_{(x, \pi) \in G_{1}(N_{h,g_{\beta}}(\gamma_{g_{\beta}}))} |JG_{\beta}|^{-1}  \\ \nonumber &+  Lip((JG_{\beta})^{-1}) \sup_{(x, \pi) \in G_{1}(N_{h,g_{\beta}}(\gamma_{g_{\beta}}))} |f| \Big) \\ \nonumber &\le  \exp(-\phi_{\beta}) \Big( (1-\beta)^{-1}  +  \sup_{N_{h,g_{\beta}}(\gamma_{g_{\beta}})} \Big| \frac{\partial \phi_{\beta}}{\partial \nu} \exp(-\phi_{\beta}) \Big| \Big) 
    \\ \nonumber &= \exp(-\phi_{\beta}) \Bigg( (1-\beta)^{-1}  +  \sup_{y \in (-h,h)} \Bigg| \frac{2ABy e^{-\frac{B}{h^{2}-y^{2}}}}{(h^{2}-y^{2})^{2}} \exp(-\phi_{\beta}) \Bigg| \Bigg)
    \\ &\le \exp(-\phi_{\beta}) ( (1-\beta)^{-1}  +  Mh \exp(-\phi_{\beta})).
\end{align}
Now consider the bound of (43) and (44). For (43), $(1-\beta)^{-1}$ converges to $1$ as $\beta$ goes to $0$. For (44), the right hand side $\exp(-\phi_{\beta}) ( (1-\beta)^{-1}  +  Mh \exp(-\phi_{\beta}))$ converges to $1+Mh$ which is smaller than $2$ as $\beta$ goes to $0$. Hence, we can take $\beta>0$ such that $\max ((1-\beta)^{-1}, \exp(-\phi_{\beta}) ( (1-\beta)^{-1}  +  Mh \exp(-\phi_{\beta})))<3$.

By this choice of $\beta$, we take some smaller tubular neighborhood $N_{h'}(\gamma_{g_{\beta}})$ such that ambient Gaussian curvature is negative and $N_{h'}(\gamma_{g_{\beta}})$ is a mean convex neighborhood for some $h'>0$ by Lemma 5.3. Note that $N_{h'}(\gamma_{g_{\beta}})$ may be smaller than $N_{h,g_{\beta}}(\gamma_{g_{\beta}})$. Let us take $\epsilon_{0} = h'^{2}/20$. By Lemma 3.1 and (27), curves $ (G_{\beta})_{\sharp}(\Phi(x)) $ are supported in $N_{h'}(\gamma_{g_{\beta}})$ for any $x \in X$. Thus, the family $\{ (G_{\beta})_{\sharp}(\Phi(x)) \}_{x \in X}$ satisfies the assumption of Lemma 6.2. From here let us remove subscript $\beta$ for the brevity of notation. By applying Lemma 6.2, we obtain
\begin{align}
    \nonumber F(&H(t,x),\gamma)\\ \nonumber &= \sup \Big\{ \Big| \int f (JG)^{-1} dH'(t,x) - \int f (JG)^{-1} dG_{\sharp}\gamma \Big| : f \in C_{c}(G_{1}(S^{2})), |f| \le 1 , Lip(f) \le 1 \Big\} \\ \nonumber &= 3 \sup \Big\{ \Big| \int \frac{f (JG)^{-1}}{3} dH'(t,x) - \int \frac{f (JG)^{-1}}{3} dG_{\sharp}\gamma \Big| : f \in C_{c}(G_{1}(S^{2})), |f| \le 1 , Lip(f) \le 1 \Big\} \\ &\le C(|\gamma|) \sqrt{\epsilon},
\end{align}
where we obtain (45) from our choice of $\beta$ above and the definition of $F$-distance (3).
\end{proof}
\end{lem}

\subsection{The smooth interpolation lemma}
We now state the interpolation lemma. Our interpolation is to construct the homotopy between two families of curves which are close to each other near the boundary of the local min-max ball with controlled length.

The idea is as follows. Suppose for each curve $\Phi(x)$ in the family $\{ \Phi(x) \}_{x \in X}$ there exists continuous $w: X \rightarrow \partial B^{j}(0,1)$ such that $F(\Phi(x),(F_{w(x)})_{\sharp}\gamma)<\epsilon$ uniformly for small $\epsilon$ where $\gamma$ is a fixed geodesic and $\{F_{v} \}_{v \in \overline{B}^{j}}$ is a local min-max diffeomorphism in Section 4. Then we consider the family $\{ (F_{w(x)})^{\sharp}(\Phi(x)) \}$ and curves in this family are sufficiently close to the geodesic $\gamma$ and we apply the $F$-distance estimate in Lemma 6.3. 

\begin{lem} Let $\gamma$ be a fixed geodesic and $X$ be a simplicial complex with finite dimension $k$. Suppose $\Phi: X \rightarrow \mathcal{S}$ be a $k$-parameter family of simple closed curves and $\{ F_{v}\}_{v \in \overline{B}^{j}}$ be a local min-max family of diffeomorphisms in Theorem 4.1. There exists $\epsilon_{1}=\epsilon_{1}(\gamma)>0$ such that satisfying the following property : If there is a continuous function $w: X \rightarrow \partial B^{j}$ satisfying \begin{equation}
F(\Phi(x),(F_{w(x)})_{\sharp}\gamma)<\epsilon
\end{equation}
with $0<\epsilon<\epsilon_{1}$, then there exists a homotopy $H: [0,1] \times X \rightarrow \mathcal{S}$ such that $H(0,x) = \Phi(x)$, $H(1,x) = (F_{w(x)})_{\sharp}\gamma$ satisfying
\begin{equation*}
    \sup \{ F(H(t,x),(F_{w(x)})_{\sharp}\gamma) : x \in X \text{ and } t \in [0,1] \} < C(|\gamma|) \sqrt{\epsilon}.
\end{equation*}
\begin{proof}
By (iii) of Theorem 4.1, let us take the local min-max ball $\{ F_{v}\}_{v \in \overline{B}^{j}}$ to satisfy
\begin{equation}
    F((F_{v})_{\sharp}(V_{1}),(F_{v})_{\sharp}(V_{2})) \le 2 F(V_{1},V_{2})
\end{equation}
for any $V_{1},V_{2} \in IV_{1}(S^{2})$. We define a pulled-back family $\Psi: X \rightarrow \mathcal{S}$ as $\Psi(x) := (F_{-w(x)})_{\sharp}(\Phi(x))$. Then the following holds by (46) and (47):
\begin{equation}
    F(\Psi(x),\gamma)  = F((F_{-w(x)})_{\sharp}(\Phi(x)), (F_{-w(x)})_{\sharp}(F_{w(x)})_{\sharp}(\gamma))\le  2F(\Phi(x), (F_{w(x)})_{\sharp}(\gamma))< 2 \epsilon. 
\end{equation}
We take $\epsilon_{1}=\epsilon_{0}/2$ where $\epsilon_{0}=\epsilon_{0}(\gamma)$ in Lemma 6.3. By Lemma 6.3 and (48), there exists a squeezing homotopy $H'(t,x)$ such that $H'(0,x) = \Psi(x)$ and $H'(1,x) = \gamma$ with
\begin{equation}
    \sup \{ F(H'(t,x),\gamma) : x \in X \text{ and } t \in [0,1] \} < C(|\gamma|) \sqrt{\epsilon}.
\end{equation}
Let us define the pushforward homotopy $H: [0,1] \times X \rightarrow \mathcal{S}$ of $H'$ as $H(t,x) = (F_{w(x)})_{\sharp}(H'(t,x))$. Then the following holds:
\begin{eqnarray}
    F(H(t,x), (F_{w(x)})_{\sharp}(\gamma)) &=& F((F_{w(x)})_{\sharp}(H'(t,x)), (F_{w(x)})_{\sharp}(\gamma)) \nonumber \\ &\le& 2F(H'(t,x),\gamma) < C(|\gamma|) \sqrt{\epsilon},
\end{eqnarray}
for any $t \in [0,1]$ and $x \in X$. The inequalities in (50) come from (47) and (49), and $H$ is a desired homotopy. Note that $C(|\gamma|)$ does not depend on $w(x)$. 
\end{proof}
\end{lem}

\section{Proof of Theorem 1.2}
In this section, we prove Theorem 1.2 on the Morse index bound. The following theorem implies there are at least three simple closed geodesics on bumpy sphere with Morse index $1,2$ and $3$. These geodesics with Morse index $1,2$ and $3$ exist generically since bumpy metrics are generic in $C^{k}$-Baire sense from \cite{Ab}. We adopt analogous notions from \cite{MN3} in our proof. We apply the interpolation lemma (Lemma 6.4) we proved in the last section (cf. Theorem 3.8 in \cite{MN3}).
\begin{thm}
Suppose $(S^{2}, g)$ is a $2$-sphere with a bumpy metric. Then for each $k=1,2,3$ there exists a simple closed geodesic $\gamma_{k}$ with
\begin{equation*}
    index(\gamma_{k}) = k
\end{equation*}
and the lengths of these three geodesics satisfy $|\gamma_{1}|<|\gamma_{2}|<|\gamma_{3}|$.
\begin{proof}
For $k=1,2,3$, $W_{L_{k}}$ is a finite set since the metric $(S^{2},g)$ is bumpy by Corollary A.2 on spheres with bumpy metric. For each geodesic $\gamma \in W_{L_{k}}$, by (iii) of Theorem 4.1, we can take the local min-max ball $\{ F_{v}\}_{v \in \overline{B}^{j}}$ to satisfy (47). Then there exists some $b(\gamma)>0$ such that  \begin{equation}
    |(F_{v})_{\sharp}(\gamma)|<L_{k}-b(\gamma)
\end{equation} for $v \in \partial B^{j}$ since the length functional is strictly concave in the local min-max ball.

We pick small $0<s<\min_{\gamma \in W_{L_{k}}} \min(\epsilon(\gamma), b(\gamma)^{2} C(L_{k})^{-2}/30 ,\epsilon_{1}(\gamma)/3)$ satisfying $B^{F}_{2s}(\gamma) \cap W_{L_{k}} = \gamma$ for every $\gamma \in W_{L_{k}}$, where $\epsilon(\gamma)$ is a constant in Theorem 4.1, $b(\gamma)$ is a constant in (51), and $\epsilon_{1}(\gamma)$ and $C(L_{k})=C(|\gamma|)$ are constants in Lemma 6.4 for $\gamma \in W_{L_{k}}$.

We can also find a neighborhood $B^{F}_{s'} (\gamma)$ of $\gamma$ such that $|(F_{v})_{\sharp}(V)|<|\gamma|-\frac{b}{2}$ for $v \in \partial B^{j}$ and $V \in B^{F}_{s'} (\gamma)$ by setting $4s'<b$:
\begin{eqnarray}
|(F_{v})_{\sharp}(V)|  &\le&  |(F_{v})_{\sharp}(\gamma)| + F((F_{v})_{\sharp}(V),(F_{v})_{\sharp}(\gamma)) \\
&\le& |(F_{v})_{\sharp}(\gamma)| + 2 F(V, \gamma) \nonumber \\
&<&  |(F_{v})_{\sharp}(\gamma)| + 2s' \nonumber\\
&<& L_{k}-b+ 2s' \nonumber\\
&<& L_{k}- \frac{b}{2} \nonumber
\end{eqnarray}
where we used the definition of $F$-metric distance (3) for (52). We replace $s$ by $s'$ if $s'<s$.

We consider a minimizing sequence $\{ \Phi_{i}  (x) \}_{x \in X}$. Here $X$ is a $k$-dimensional simplicial complex. We apply Theorem 2.3 to obtain a tightened sequence $\{ \hat{\Phi}_{i} (x)\}_{x \in X}$ such that
\begin{equation}
        \{ \hat{\Phi}_{i}(x) \in IV_{1}(S^{2}) : |\hat{\Phi}_{i}(x)| \ge L_{k}-a \} \in \bigcup_{\gamma \in \Lambda(\{ \Phi_{i} \}) \cap W_{L_{k},k}} B^{F}_{s}(\gamma)
        \end{equation}
for some $0<a<L_{k}$ and for all sufficiently large $i$. Note that the critical set $\Lambda(\{ \hat{\Phi}_{i} \})$ consists of multiplicity one geodesics and does not contain any geodesic whose Morse index is larger than or equal to $k+1$.

Let $\Lambda (\{ \hat{\Phi}_{i} \})= \{ \tilde{\gamma}_{1}, ..., \tilde{\gamma}_{q} \}$. We prove that there exists an element $\tilde{\gamma} \in \Lambda_{k} (\{ \hat{\Phi}_{i} \})$ with $index(\tilde{\gamma}) = k$. We show this by contradiction, thus let us assume that $index(\tilde{\gamma}_{i})< k$ for every $i=1, ..., q$. We consider a $k$-sweepout restricted on the union of $k$-faces $\tilde{X}^{k}$ of $k$-skeleton $X^{k}$ of $X$ as in \cite{MN3}.

Since $\{ \hat{\Phi}_{i} \}$ is a sequence of $k$-sweepout, we can let there is a $\sigma \in H_{k}(\tilde{X}^{k}, \mathbb{Z}_{2})$ such that
\begin{equation*}
    \hat{\Phi}_{i}^{*} (\overline{\lambda}^{k})|_{\tilde{X}^{k}} \cdot \sigma = 1,
\end{equation*}
where $0 \neq \overline{\lambda} \in H^{1}(\mathcal{S}, \mathbb{Z}_{2})$. Also we can set $\sigma = [\sum^{h}_{l=1} t_{l} ]$ for some $k$-dimensional simplices $t_{1}, ..., t_{h}$ in $\tilde{X}^{k}$ from the equivalence between simplicial homology and singular homology where $\sum_{l=1}^{h} \partial t_{l}=0$ holds.

Let $Y_{i}$ be a $m_{i}$'th successive barycentric subdivision of $\cup^{h}_{l=1} t_{l}$ so that $F(\hat{\Phi}_{i}(x),\hat{\Phi}_{i}(y))<a/2$ whenever $x,y$ is in a same simplex in $\cup^{h}_{l=1} t_{l}$. Denote $W_{i}$ be a union of all $k$-dimensional simplices $t \in Y_{i}$ such that $|\hat{\Phi}_{i}(x)| \ge L_{k}- a/2$ for some $x \in t$. Then $|\hat{\Phi}_{i}(y)| \ge L_{k}- a$ for every $y \in W_{i}$. Let $W_{i,1}, ...,W_{i,r}$ be the connected components of $W_{i}$. Then from (53), for each $1 \le p \le r$ there exists $1 \le q_{p} \le q$ with
\begin{equation}
    F(\hat{\Phi}_{i}(y), \gamma_{q_p}) <s
\end{equation}
for every $y \in W_{i,p}$. If $\partial (\sum_{t \in W_{i}}t) \neq 0$, the following holds:
\begin{equation}
    L_{k}-a \le |\hat{\Phi}_{i}(y)| \le L_{k}-a/2
\end{equation}
for $y \in \partial (\sum_{t \in W_{i}}t)$ since $y$ belongs to both $W_{i}$ and $Y_{i} \setminus W_{i}$.

Now we fix one connected component $W_{i,p}$. For $y \in supp(\partial (\sum_{t \in W_{i,p}}t))$, we define the function $L^{y} : \overline{B}^{j} \rightarrow [0,\infty)$ by
\begin{equation*}
    L^{y}(v) = |(F_{v})_{\sharp} (\hat{\Phi}_{i}(y))|.
\end{equation*}
By our choice of $s$ and Theorem 4.1, the function $L^{y}$ is strictly concave and has a unique maximum at $m_{i}(y) \in B^{j}_{\frac{1}{2}} (0)$. The function $m_{i} : supp(\partial(\sum_{t \in W_{i,p}}t)) \rightarrow B^{j}_{\frac{1}{2}} (0)$ is moreover continuous. By applying Theorem 4.1 with $V= \hat{\Phi}_{i}(y)$, we have $m_{i}(y) \neq 0$ for every $y \in supp(\partial (\sum_{t \in W_{i,p}}t))$ from the upper bound of (55).

First we consider the case of $\partial(W_{i,p})=0$. We apply Lemma 5.4 to squeeze the family of curves in the small ball. More precisely,
\begin{eqnarray}
    1 &=& \overline{\lambda}^{k} \cdot \Bigg[ \sum_{t \in Y_{i}} (\hat{\Phi}_{i}) _{\sharp} (t) \Bigg] \nonumber \\
    &=& \overline{\lambda}^{k} \cdot \Bigg[ \sum_{t \in Y_{i} \setminus W_{i,p}} (\hat{\Phi}_{i}) _{\sharp} (t) \Bigg] + \overline{\lambda}^{k} \cdot \Bigg[ \sum_{t \in W_{i,p}} (\hat{\Phi}_{i}) _{\sharp} (t) \Bigg] \\
    &=& \overline{\lambda}^{k} \cdot \Bigg[ \sum_{t \in Y_{i} \setminus W_{i,p}} (\hat{\Phi}_{i}) _{\sharp} (t) \Bigg]
\end{eqnarray}
(56) comes from $\partial(W_{i,p})=0$, and (57) is from (54) and Lemma 5.4.

If $j=0$, $|\hat{\Phi}_{i}(y)| \ge L_{k}$ holds by Theorem 4.1 for every $y \in W_{i,p}$. Suppose $\partial(W_{i,p}) \neq 0$ then $|\hat{\Phi}_{i}(y)| \le L_{k}- a/2$ for $y \in \partial(W_{i,p})$ by (51). However, this contradicts $|\hat{\Phi}_{i}(y)| \ge L_{k}$ for $y \in W_{i,p}$ and we conclude $\partial(W_{i,p})=0$.

Suppose $\partial (\sum_{t \in W_{i,p}}t) \neq 0$ and $j \ge 1$. Recall that $m_{i}(y) \neq 0$ for $y \in supp(\partial(\sum_{t \in W_{i,p}}t))$. We apply a one-parameter length-decreasing flow $\{ \phi^{y}(\cdot,t) \} _{t \ge 0} \subset Diff(\overline{B}^{j})$ generated by
\begin{equation}
    v \mapsto -(1-|v|^{2}) \nabla L^{y}(v)
\end{equation}
to the curves $\hat{\Phi}_{i}(supp(\partial(\sum_{t \in W_{i,p}}t)))$. Note that $\lim _{t \rightarrow \infty} \phi^{y}(v,t)\in \partial B^{j}$ and the limit is uniform for each connected component of $supp(\partial(\sum_{t \in W_{i,p}}t))$ if $|v- m_{i}(y)| \ge \eta$ for some $\eta >0$.

We construct a first homotopy
\begin{equation}
    H_{1} : [0,1] \times  supp \Big(\partial(\sum_{t \in W_{i,p}}t )\Big) \rightarrow IV_{1}(S^{2} )
\end{equation}
that for some large $t_{0}>0$, which is defined by
\begin{equation*}
    H_{1}(t,y) = (F_{\phi^{y}(0,t_{0}t)})_{\sharp} (\hat{\Phi}_{i}(y))
\end{equation*}
satisfying
\begin{enumerate}
    \item $H_{1}(0,y) = \hat{\Phi}_{i}(y)$,
    \item $F(H_{1}(1,y), (F_{w(y)})_{\sharp}(\gamma_{q_{p}}))< 3s$ for some continuous function 
    \begin{equation*}
        w: supp \Big(\partial(\sum_{t \in W_{i,p}}t) \Big) \rightarrow \partial B^{j}
    \end{equation*}
    \item $|H_{1}(t,y)|< L_{k} - a/2$ 
\end{enumerate}
for $y \in supp (\partial(\sum_{t \in W_{i,p}}t))$ and $t \in [0,1]$. Then for the condition (2),
\begin{eqnarray}
F(H_{1}(1,y), (F_{w(y)})_{\sharp}(\gamma_{q_{p}})) &=& F((F_{\phi^{y}(0,t_{0})})_{\sharp} (\hat{\Phi}_{i}(y)),(F_{w(y)})_{\sharp}(\gamma_{q_{p}})) \nonumber\\
&\le& F((F_{\phi^{y}(0,t_{0})})_{\sharp} (\hat{\Phi}_{i}(y)),(F_{w(y)})_{\sharp}(\hat{\Phi}_{i}(y)))\nonumber \\ &+& F((F_{w(y)})_{\sharp}(\hat{\Phi}_{i}(y)),(F_{w(y)})_{\sharp}(\gamma_{q_{p}})) \nonumber \\
&\le& s + 2F(\hat{\Phi}_{i}(y),\gamma_{q_{p}}) < 3s
\end{eqnarray}
for large $t_{0}$ and (60) holds since $\lim _{t \rightarrow \infty} \phi^{y}(v,t)\in \partial B^{j}$ and the convergence is uniform. By the fact that $\{ \phi^{y}(\cdot,t) \}$ induces a length-decreasing flow and (55), the condition (3) holds.

We now apply Lemma 6.4 to construct a homotopy between $\{ H_{1}(1,y) \}$ and $\{ (F_{w(y)})_{\sharp}(\gamma_{q_{p}})\}$. Let us put $\epsilon = 3s<\epsilon_{1}$ in Lemma 6.4. Then by Lemma 6.4, there exists a homotopy
\begin{equation}
    H_{2} : [0,1] \times  supp \Big(\partial(\sum_{t \in W_{i,p}}t )\Big) \rightarrow IV_{1}(S^{2} )
\end{equation}
such that
\begin{enumerate}
    \item $H_{2}(0,y) = H_{1}(1,y)$,
    \item $H_{2}(1,y) = (F_{w(y)})_{\sharp}(\gamma_{q_{p}})$
    \item $F(H_{2}(t,y),(F_{w(y)})_{\sharp}(\gamma_{q_{p}}))< b/3$ and so $|H_{2}(t,y)|< L_{k}- b/6$ 
\end{enumerate}
for $y \in supp (\partial(\sum_{t \in W_{i,p}}t))$ and $t \in [0,1]$. The $F$-distance estimate of the condition (3) comes from
\begin{equation}
F(H_{2}(t,y),(F_{w(y)})_{\sharp}(\gamma_{q_{p}})) \le C(L_{k}) \sqrt{3s} \le C(L_{k}) \sqrt{\frac{3b^{2}}{30C(L_{k})^{2}}} < \frac{b}{3},
\end{equation}
where $C(L_{k})$ is the constant from Lemma 6.4. The first inequality of (62) is the conclusion of Lemma 6.4 and the second inequality is from the choice of $s$. The length bound directly comes from the length estimate (52) on the boundary of local min-max ball and the definition of $F$-distance (3):
\begin{equation*}
    |H_{2}(t,y)| \le |(F_{w(y)})_{\sharp}(\gamma_{q_{p}})| + F(H_{2}(t,y),(F_{w(y)})_{\sharp}(\gamma_{q_{p}})) < L_{k}- \frac{b}{2}+ \frac{b}{3} = L_{k}- \frac{b}{6}.
\end{equation*}

Since $H_{k-1}(\partial B^{j}, \mathbb{Z}_{2}) = 0$ for $j<k$, we have 
\begin{equation*}
    \Bigg[w_{\sharp} \Big(\partial (\sum_{t \in W_{i,p}}t) \Big)\Bigg] = 0.
\end{equation*}
Hence there is a $k$-dimensional singular chain $\sum_{j} \alpha_{j}$ on $\partial B^{j}(0,1)$ such that $\sum_{j} \partial \alpha_{j} = w_{\sharp}(\partial (\sum_{t \in W_{i,p}}t))$. Then we have the associated singular $k$-simplex $\hat{\alpha}_{j}: \Delta^{k} \rightarrow IV_{1}(S^{2})$ as
\begin{equation*}
    \hat{\alpha}_{j}(y) = F_{\alpha_{j}(y)}( \gamma_{q_{p}}),
\end{equation*}
where $y \in \Delta^{k}$.

Now we define a new singular chain. We denote $W_{i}^{0}$ as a union of all components $W_{i,p}$ with $\partial (\sum_{t \in W_{i,p}}t) = 0$. Then we consider the singular chain
\begin{equation*}
    z_{i,p} = \sum _{t \in W_{i,p}} (\hat{\Phi}_{i})_{\sharp} (t) + (H_{1}+H_{2})_{\sharp} \Bigg([0,1] \times supp \Big(\partial(\sum_{t \in W_{i,p}}t) \Big) \Bigg) + \sum_{j} \hat{\alpha}_{j}.
\end{equation*}
Then $\partial z_{i,p} =0$. Moreover, following holds by Lemma 5.4 and Remark 5.5 since our perturbation were along some small tubular neighborhood $N_{\alpha}(\gamma_{q_{p}})$ of $\gamma_{q_{p}}$ and can be squeezed to the center circle $\gamma_{q_{p}}$ so that $z_{i,p}$ is homologically trivial:
\begin{equation*}
    \overline{\lambda}^{k} \cdot [z_{i,p}] =0. 
\end{equation*}
Hence, the $k$-dimensional cycle
\begin{equation}
  \tilde{z} = \sum _{t \in Y_{i} \setminus W_{i}^{0} } (\hat{\Phi}_{i})_{\sharp}(t) + \sum_{(i,p) \in I} z_{i,p}, 
\end{equation}
where $I := \{ (i,p): index(\gamma_{q_{p}})>0 \,\, and \,\, \partial(W_{i,p}) \neq 0\}$ satisfies $\overline{\lambda}^{k} \cdot [\tilde{z}] =1$ and so (63) is a $k$-sweepout having $|W| \le L_{k}- \min(b/6,a/2) $ for every $W \in image(\tilde{z})$. Then we mollify $\tilde{z}$ to make deformations depend smoothly of the parameters and we obtain a smooth $k$-sweepout $\tilde{\tilde{z}}$ with $|\tilde{W}| \le L_{k}- \min(b/6,a/2) $ for every $\tilde{W} \in image(\tilde{\tilde{z}})$. Now we can construct a $\Delta$-complex $\tilde{\tilde{Z}}$ from $\tilde{\tilde{z}}$ and a continuous map $\Xi_{i} : \tilde{\tilde{Z}} \rightarrow \mathcal{S}$ in the smooth topology such that $(\Xi_{i})^{*}(\overline{\lambda}^{k}) \neq 0$ and $|\Xi_{i}(y)| \le L_{k} - \min(b/6,a/2)$ for any $y \in \tilde{\tilde{Z}}$. This gives contradiction to the fact that the critical length is $L_{k}$. The distinction among length $L_{k}$ by Corollary 2.2 gives the latter conclusion.

\end{proof}
\end{thm}
\begin{rmk} Since we chose $a$ and $s$ after taking an appropriate local min-max diffeomorphism $\{F_{v} \}_{v \in \overline{B}^{k}}$ and $b$ for length bound, the length along our interpolation could be bounded with the number which is strictly smaller than the width.
\end{rmk}
We now prove Morse index characterization in general spheres (Corollary 1.4) as following:

\begin{cor}
For a 2-sphere $(S^{2},g)$ with a smooth metric, for $k=1,2,3$ there exists a closed and embedded geodesic $\gamma_{k}$ with
\begin{equation*}
    index(\gamma_{k}) \le k \le index(\gamma_{k})+ nullity (\gamma_{k}).
\end{equation*}
\begin{proof}
We can approximate $(S^{2},g)$ by a sequence of bumpy metrics $\{ (S^{2},g_{i}) \}$ in the $C^{\infty}$-sense by generic property of bumpy metric by \cite{Ab}. Let $\gamma_{k,i}$ be an element of the set of simple closed geodesics obtained by the theorem above from $k$-sweepouts on $(S^{2},g_{i})$. We obtain a geodesic $\gamma_{k}$ in $(S^{2},g)$ as a subsequential limit of $\gamma_{k,i}$ by the local compactness theorem and local uniqueness theorem of solution of ODE. The local convergence and the compactness give convergence toward the geodesic $\gamma_{k}$.

Moreover, we see the convergence of eigenvalues through the convergence of geodesics by the variational characterization of eigenvalues of the stability operator. The stability operator of $\gamma_{k,i}$ with normal variation along $f:\gamma_{k,i} \rightarrow \mathbb{R}$ is (with the normal coordinate system along geodesics)
\begin{equation*}
    \delta_{\gamma_{k,i}}(f) = \int_{\gamma_{k,i}} (|\nabla_{\gamma_{k,i}}f|^{2} - K f^{2}) ds.
\end{equation*}
By the $C^{\infty}$-convergence of geodesics and Gaussian curvature, and a variational characterization of $j$-th eigenvalues in terms of Rayleigh quotients as
\begin{equation*}
    \lambda_{k,i,j} = \inf_{j\text{-plane}\,\, P \subset W^{1,2}(\gamma_{k,i})\setminus \{ 0\} } \max_{f \in P} \Bigg\{ \frac{\delta_{\gamma_{k,i}}(f)}{\int f^{2}}\Bigg\}
\end{equation*}
(cf. Lemma 1.34 of \cite{CM}), we have the convergence of $j$-th eigenvalues $\lambda_{k,i,j}$ of stability operators of $\gamma_{k,i}$ to $j$-th eigenvalue $\lambda_{k,j}$ of stability operators of $\gamma_{k}$ as $i$ goes to infinity. Then we obtain $\lambda_{k,k} \le 0$ since $\lambda_{k,i,k} < 0$ for $k=1,2,3$ by the assumption. This gives the conclusion that $index(\gamma_{k}) \le k \le index(\gamma_{k})+ nullity (\gamma_{k})$ so that we obtained desired geodesics.  
\end{proof}
\end{cor}

\section{Morse Inequalities}

The goal of this section is to obtain Morse inequalities for the length functional of simple closed geodesics with fixed homotopy class on orientable surfaces. We adopt arguments in \cite{MMN} and explain differences mainly.

Let $\Pi$ be a fixed homotopy class of simple closed curves on $M$. We can deduce a strong Morse inequality for a fixed homotopy class $\Pi$ and so on the whole space of simple closed curves with bounded length. Also we denote $\mathcal{S}$ as the set of (nonparametrized) embedded smooth curves on $M$ in this section. We define $\mathfrak{M}_{g,\Pi}$ as a collection of (multiplicity one) simple closed geodesics with fixed homotopy class $\Pi$ in $(M^{2},g)$. In this section, let us suppose that the Riemannian metric $g$ is bumpy in $M^{2}$ which is $C^{\infty}$-generic.

For $k \in \mathbb{Z}_{+} $ and $a \in (0,\infty)$, let $c_{k}(a,\Pi)$ be the number of simple closed geodesics with $index(\gamma)=k$, $length(\gamma)<a$ and $\gamma \in \mathfrak{M}_{g,\Pi}$. Let also $b_{k}(a,\Pi)$ be a $k$-th Betti number of the space of embedded curves
\begin{equation*}
    \mathcal{Z}^{a}_{\Pi} = \{ \gamma \in \Pi : |\gamma| < a \}
\end{equation*}
with the varifold $F$-metric. Note that $\mathcal{Z}^{a}_{\Pi}$ is open in $\Pi$ for any $a$. Also denote that $b_{k}(a)$ and $c_{k}(a)$ to be $k$-th Betti number of the space of embedded curves $\mathcal{Z}^{a} = \{ \gamma \in \mathcal{S} : |\gamma| < a \}$ and a number of simple closed geodesics with $index(\gamma)=k$, $length(\gamma)<a$, respectively.

We start with the following lemma by the compactness theorem (Theorem A.1) and the finiteness of simple closed geodesics with bounded length (Corollary A.2):
\begin{lem}
$\sum_{\Pi} \sum_{k=0}^{\infty} c_{k}(a,\Pi) <\infty$ so that $\sum_{k=0}^{\infty} c_{k}(a,\Pi) <\infty $ for every $a>0$.
\end{lem}
We have the following analogue of Proposition 3.3 directly from \cite{MMN}.
\begin{prop}
For a $C^{\infty}$-generic Riemannian metric $g$ on $M$, we have
\begin{itemize}
    \item every simple closed geodesic is nondegenerate;
    \item and if
    \begin{equation*}
        p_{1} \cdot length_{g}(\gamma_{1}) + \hdots + p_{N} \cdot length_{g}(\gamma_{N}) = 0,
    \end{equation*}
    with $\{ p_{1},\hdots, p_{N} \} \subset \mathbb{N}$ and $\{ \gamma_{1}, \hdots, \gamma_{N} \} \subset \mathfrak{M}_{g,\Pi}$, and $\gamma_{k} \neq \gamma_{l}$ whenever $k \neq l$, then
    \begin{equation*}
        p_{1} = \hdots = p_{N} =0.
    \end{equation*}
\end{itemize}
\end{prop}

We now see the relation between the relative homology $H_{k}(\mathcal{Z}^{t}_{\Pi},\mathcal{Z}^{s}_{\Pi})$ with $\mathbb{Z}_{2}$-coefficients and numbers of simple closed geodesics. To obtain an analogue of Homology Min-Max Theorem 3.5 of \cite{MMN}, we adopt the arguments to obtain the lower bound of Morse index of geodesic in the proof of Theorem 7.1, and apply curve shortening flow in place of area decreasing homotopies. Note that Local min-max theorem (Theorem 4.1) and the interpolation lemma (Lemma 6.4) still can be applied even for general geodesics.  

\begin{thm}[Homology Min-max Theorem] Let $\sigma \in H_{k}(\mathcal{Z}_{\Pi}^{t},\mathcal{Z}_{\Pi}^{s})$ be a nontrivial homology class, $0 \le s < t$, and let
\begin{equation*}
    W(\sigma) = \inf_{[\sum s_{i}] = \sigma} \sup_{i,x \in \Delta^{k}} |s_{i}(x)|.
\end{equation*}
Suppose that, for some $\epsilon>0$, there is no $\gamma' \in \mathfrak{M}_{g,\Pi}$ with $length(\gamma') \in (s- \epsilon,s)$ and $index(\gamma') \le k-1$. Then $W(\sigma) \in [s,t)$ with $W(\sigma) >0$ if $s=0$. Moreover there exists $\gamma \in \mathfrak{M}_{g,\Pi}$ with $index(\gamma) =k$ and
\begin{equation*}
    length(\gamma) = W(\sigma).
\end{equation*}

\end{thm}
We obtain following proposition which is the analogue of Proposition 3.6 in \cite{MMN}. An inspection of the proof gives that we can prove the following proposition by replacing the interpolation results in Almgren-Pitts setting by smooth interpolation lemmas (Lemma 6.3 and Lemma 6.4).  
\begin{prop}
Suppose that the metric $g$ satisfies the conditions in Proposition 8.2. Let $\gamma \in \mathfrak{M}_{g,\Pi} \cap \mathcal{S}$ with $index(\gamma) =k$ and $length(\gamma)=l$. For every $\epsilon>0$, there exists $l-\epsilon<s<l<t<l+\epsilon$ such that
\begin{eqnarray}
b_{i}(\mathcal{Z}^{t}_{\Pi},\mathcal{Z}^{s}_{\Pi}) &=& 0 \text{  if $i \neq k$} \\
b_{k}(\mathcal{Z}^{t}_{\Pi},\mathcal{Z}^{s}_{\Pi}) &=& 1.
\end{eqnarray}
\begin{proof}
We describe the outline of the proof of Proposition 3.6 of \cite{MMN} and mainly discuss about the modifications one needs to make. 

There is a $\delta>0$ such that $\gamma$ is the unique element of $\mathfrak{M}_{g,\Pi}$ whose length is in $(l-\delta,l+\delta)$ since our choice of the metric satisfies Proposition 8.2. (64) follows from Theorem 8.3 directly.

Let us consider a local min-max family associated with $\gamma$. We take the functional $L^{*}$ in (11). As in the argument of the proof of Theorem 4.1, we can pick $\epsilon_{2}>0$ such that for any $V \in IV_{1}(S^{2})$ induced by simple closed curves with $F(V,\gamma) < \epsilon_{2}$, we have $L^{*}(V) > L^{*}(\gamma)$ unless $\gamma = V$. Since $L^{*}$ is strictly convex by Theorem 5 in \cite{W2}, we can pick $\kappa = \kappa(\alpha, \gamma, \epsilon_{2})>0$ such that for $V$ induced by simple closed curves with $F(V,\gamma) \le \epsilon_{2}$ and $L^{*}(V) \le L^{*}(\gamma) + \kappa$, then $F(V,\gamma) \le \alpha$. Let us choose $t$ to be close to $l$ such that $t-l< \min \{\delta, \kappa(\epsilon_{2}/16, \gamma, \epsilon_{2}), 2 \sqrt{\epsilon_{0}}C(l) \}$, where $\epsilon_{0}$ and $C(l)$ are constants from Lemma 6.3. Then we take a local min-max family $\{ F_{v} \}_{v \in \overline{B}^{k}}$ satisfying (47) and 
\begin{equation}
    F((F_{v})_{\sharp}(\gamma), \gamma) \le \min (\epsilon_{2}/4, (t-l)^{2}/16C'(l)^{2}),
\end{equation}
where we take $C'(l)$ as a larger one between two constants $C(l)$ from Lemma 6.3 and Lemma 6.4 and for all $v \in \overline{B}^{k}$.

Moreover, let us choose $0<b< \min(\delta, (t-l)^{2}/4C(l)^{2})$ such that
\begin{equation}
|(F_{v})_{\sharp}(\gamma)| <l - b
\end{equation}
for any $v \in \partial B^{k}$.

Suppose that $\overline{\Phi}:\overline{B}^{k} \rightarrow IV_{1}(S^{2})$ is defined as $\overline{\Phi}(v) = (F_{v})_{\sharp}(\gamma)$, then we define the relative homology class as following:
\begin{equation*}
    \overline{\sigma} = \overline{\Phi}_{*}([\overline{B}^{k}]) \in H_{k}(\mathcal{Z}^{t},\mathcal{Z}^{l-b/4}),
\end{equation*}
where $[\overline{B}^{k}]$ is a generator of $H_{k}(\overline{B}^{k}, \partial \overline{B}^{k})$. We obtain the following claim by following the proof of Claim 3.7 in \cite{MMN}.
\begin{claim}
$\overline{\sigma} \neq 0.$
\end{claim}
The proposition directly follows from the following claim:
\begin{claim}
    $H_{k}(\mathcal{Z}^{t},\mathcal{Z}^{l-b/4}) = \{0, \overline{\sigma} \}.$
    \begin{proof}
        Let $0 \neq \sigma \in H_{k}(\mathcal{Z}^{t},\mathcal{Z}^{l-b/4})$. By Theorem 8.3, $W(\sigma) = l$. We can suppose that there is a minimizing sequence $\{ \Phi_{j} \}$ of maps $\Phi_{j} : X \rightarrow \mathcal{Z}^{t}$ such that $\Phi_{j}(\partial X) \subset \mathcal{Z}^{l-b/4}$, where $X$ is a $k$-dimensional $\Delta$-complex and whose width is $W(\sigma) = l$ and $(\Phi_{j})_{*}([X]) = \sigma$. By Proposition 8.2, the only element of the critical set $\Lambda(\{ \Phi_{j} \})$ is $\gamma$. By applying Theorem 2.3, we obtain a pulled-tight minimizing sequence $\{ \hat{\Phi}_{j}(x)\}$ such that there exists $\eta>0$ for given small $s>0$ satisfying
        \begin{equation}
            \{ \hat{\Phi}_{j}(x) \in IV_{1}(S^{2}) : |\hat{\Phi}_{j}(x)| \ge l-\eta \} \subset  B^{F}_{s}(\gamma)
        \end{equation}
        for sufficient large $i$ and let us take $0<s< \min(b^{2}C'(l)^{-2}/16, (t-l)^{2}C'(l)^{-2}/16,\epsilon_{1}/3)$ and $0<\eta <b/8$, where $\epsilon_{1}$ is a constant from Lemma 6.4.
        
        We fix $i$ and take a barycentric subdivision of $X$ such that $F(\Phi_{i}(x), \Phi_{i}(y))< \eta/2$ whenever $x$ and $y$ are in the same simplex of the subdivision. We define $V$ as the union of all $k$-simplicies $t \in X$ such that there exists $x \in t$ with $|\hat{\Phi}_{i}(x)| \ge l- \eta/2$. For any $x' \in t \in V$, we have $|\hat{\Phi}_{i}(x')| \ge l - \eta$. By our choice of $b$ and $\eta$, $V \cap \partial X = \empty$. In particular, we have the following length bound for every $y \in \partial V$:
        \begin{equation*}
            l- \eta \le |\hat{\Phi}_{i}(y)| \le l- \eta/2.
        \end{equation*}
        By modifying homotopies in (59) and (61) (after changing constants), we have the homotopy $H_{1} : [0,1] \times \partial V \rightarrow IV_{1}(S^{2})$ induced by length-decreasing flow $\{ \phi^{y}(\cdot,t) \} _{t \ge 0} \subset Diff(\overline{B}^{k})$ in (58) satisfying
        \begin{enumerate}
            \item $H_{1}(0,y) = \hat{\Phi}_{i}(y)$
            \item $F(H_{1}(1,y), (F_{w(y)})_{\sharp}(\gamma)) < \eta/4$ for some continuous function $w: \partial V \rightarrow \partial B^{k}$
            \item $|H_{1}(t,y)| < l - \eta/2$
        \end{enumerate}
        for every $(t,y) \in [0,1] \times \partial V$. Also, we have following $F$-distance estimate 
        \begin{eqnarray}
F(H_{1}(t,y), \gamma) &=& F((F_{\phi^{y}(0,t_{0}t)})_{\sharp} (\hat{\Phi}_{i}(y)),\gamma) \nonumber\\
&\le& F((F_{\phi^{y}(0,t_{0}t)})_{\sharp} (\hat{\Phi}_{i}(y)),(F_{\phi^{y}(0,t_{0}t)})_{\sharp} (\gamma))\nonumber + F((F_{\phi^{y}(0,t_{0}t)})_{\sharp} (\gamma),\gamma) \nonumber \\
&\le&  2s + (t-l)^{2}/16C'(l)^{2} \le (t-l)^{2}/4C(l)^{2}.
\end{eqnarray}
        Then we have a homotopy $H_{2} : [0,1] \times \partial V \rightarrow IV_{1}(S^{2})$ such that
        \begin{enumerate}
            \item $H_{2}(0,y) = H_{1}(1,y)$ and $H_{2}(1,y) = (F_{w(y)})_{\sharp}(\hat{\Phi}_{i}(y))$
            \item $F(H_{2}(t,y), (F_{w(y)})_{\sharp}(\gamma)) < b/2$ so $|H_{2}(t,y)|< l- b/2$
        \end{enumerate}
        for every $(t,y) \in [0,1] \times \partial V$. In particular, we have
\begin{eqnarray}
F(H_{2}(t,y), \gamma) &\le& F(H_{2}(t,y), (F_{w(y)})_{\sharp}(\gamma)) + F((F_{w(y)})_{\sharp}(\gamma),\gamma) \nonumber \\
&\le&  b/2 + (t-l)^{2}/16C(l)^{2} \le (t-l)^{2}/4C(l)^{2}.
\end{eqnarray}
        We define $Q$ as a cone over $\partial V$ and $\hat{w} : Q \rightarrow \overline{B}^{k}$ as a continuous map sending the vertex of the cone to origin and such that $\hat{w}(y)=w(y)$ for every $ y \in \partial V = \partial Q$.

        We consider the following $\Delta$-complex $C$ without boundary given by
        \begin{equation*}
            C:= V \cup ([0,2] \times \partial V ) \cup Q.
        \end{equation*}
        We define the continuous map $\Psi : C \rightarrow IV_{1}(S^{2})$ as
        \begin{equation*}
            \Psi(x) = \hat{\Phi}_{i}(x),
        \end{equation*}
        for $x \in V$,
        \begin{equation*}
            \Psi(s,y) = H_{1}(s,y)
        \end{equation*}
        for $(s,y) \in [0,1] \times \partial V$,
        \begin{equation*}
            \Psi(s,y) = H_{2}(s,y)
        \end{equation*}
        for $(s,y) \in [1,2] \times \partial V$, and 
        \begin{equation*}
            \Psi(q) = (F_{\hat{w}(q)})_{\sharp}(\gamma)
        \end{equation*}
        for $q \in Q$.
        We now construct a homotopy $H_{3}$ between $\Psi(C)$ and $\gamma$ in $\mathcal{Z}^{t}$ to show that
        \begin{equation}
            \Psi_{*}([C]) = 0 \in H_{k}(\mathcal{Z}^{t}, \mathcal{Z}^{l-b/2}).
        \end{equation}
        By (66), (68), (69), (70), we can apply Lemma 6.3 for $\epsilon = (t-l)^{2}/4C(l)^{2}$, and obtain the homotopy $H_{3}: [0,1] \times C \rightarrow IV_{1}(S^{2})$ such that
        \begin{enumerate}
            \item $H_{3}(0,p) = \Psi(p)$ and $H_{3}(1,p) = \gamma$
            \item $F(H_{3}(t,p), \gamma) < (t-l)/2$
        \end{enumerate}
        for every $(t,p) \in [0,1] \times C$. The definition of $F$-distance (3) gives
        \begin{equation*}
            |H_{3}(t,p)| \le |\gamma| +F(H_{3}(t,p),\gamma) \le l + (t-l)/2 < t 
        \end{equation*}
        for every $(t,p) \in [0,1] \times C$ and we have proven (71). We can finish the proof by following the remaining parts of proof of Claim 3.8 in \cite{MMN}.
        \end{proof}
\end{claim}
\end{proof}
\end{prop}
By combining Propositions and Theorems above, we obtain the strong Morse inequality in a fixed homotopy class $\Pi$.

\begin{thm} 
For each $a \in (0,\infty)$, $b_{k}(a,\Pi) < \infty$ for every $k \in \mathbb{Z}_{+}$ and the Strong Morse inequality for simple closed geodesics in a fixed homotopy class $\Pi$ hold:
\begin{equation*}
    c_{r}(a, \Pi) - c_{r-1}(a, \Pi) +\hdots+ (-1)^{r}c_{0}(a,\Pi) \ge b_{r}(a,\Pi) - b_{r-1}(a,\Pi) +...+ (-1)^{r}b_{0}(a,\Pi) 
\end{equation*}
for every $r \in \mathbb{Z}_{+}$. In particular,
\begin{equation*}
    c_{r}(a,\Pi) \ge b_{r}(a,\Pi)
\end{equation*}
for every $r \in \mathbb{Z}_{+}$.
\end{thm}
\begin{rmk} Let us denote $\mathcal{S}_{\Pi}$ to be the space of embedded closed curves in a fixed homotopy class $\Pi$. Due to the classification results of embedding of $S^{1}$ to closed surfaces (See Gramain \cite{Gra} and also follows from Grayson \cite{Gr}), if $\Pi$ contains an essential curve, then $\mathcal{S}_{\Pi}$ retracts to $S^{1}$ if $M$ is a torus, and is contractible if $genus(M) \ge 2$. In case of $\Pi$ contains a trivial curve, $\mathcal{S}_{\Pi}$ is homotopy equivalent to the unit tangent bundle to $M$, for instance if $M=S^{2}$ then $\mathcal{S}_{\Pi}$ is homotopy equivalent to $\mathbb{R}P^{3}$. This implies $b_{k}(a,\Pi)$ approaches to Betti numbers of $S^{1}$, trivial space, and unit tangent bundle of $M$ with respect to $M$ and $\Pi$ as above, as $a$ goes to infinity.   
\end{rmk}

By adding up the terms for all homotopy class by finiteness, we obtain
\begin{cor} For each $a \in (0,\infty)$, $b_{k}(a) < \infty$ for every $k \in \mathbb{Z}_{+}$ and the Strong Morse inequality for simple closed geodesics hold:
\begin{equation*}
    c_{r}(a) - c_{r-1}(a) +\hdots+ (-1)^{r}c_{0}(a) \ge b_{r}(a) - b_{r-1}(a) +...+ (-1)^{r}b_{0}(a) 
\end{equation*}
for every $r \in \mathbb{Z}_{+}$. In particular,
\begin{equation*}
    c_{r}(a) \ge b_{r}(a)
\end{equation*}
for every $r \in \mathbb{Z}_{+}$.
\end{cor}
\begin{rmk}
We restrict our surfaces by orientable surfaces due to the technical reason. To obtain a quantitative estimate as Lemma 6.4, we relied on the construction of squeezing homotopy in Lemma 5.4. In this procedure, we used the fact that a tubular neighborhood of simple closed geodesic is a cylinder, which guarantees the existence of a mean convex neighborhood of a strictly stable simple closed geodesic after a suitable deformation.
\end{rmk}
\appendix \section{Compactness theorem of closed geodesic with bounded length}
In this section, we introduce a one-dimensional analog of Sharp's compactness theorem \cite{Sh} without an index condition.

\begin{thm}
Let $M$ be a closed surface. If $\{ \gamma_{k} \}$ is a sequence of closed, connected and embedded geodesics with
\begin{equation*}
    \mathcal{H}^{1}(\gamma_{k}) \le L < \infty
\end{equation*}
for some fixed constant $L \in \mathbb{R}$ independent of $k$. Then up to subsequence, there exists a closed, connected and embedded closed geodesic $\gamma$ where $\gamma_{k} \rightarrow \gamma$ in the varifold sense with
\begin{equation*}
    \mathcal{H}^{1}(\gamma) \le L < \infty
\end{equation*}
and convergence is smooth and graphical for all $x \in M$. The multiplicity of convergence $m$ is 1 if $M$ is an orientable surface and is at most $2$ if $M$ is an unorientable surface. Moreover, if $\gamma_{k} \cap \gamma = \emptyset$ eventually, then $\gamma$ is stable, and $index (\gamma) \ge 1$ otherwise.
\begin{proof}
From Allard's compactness theorem \cite{A}, we know the existence of $\gamma$ such that (up to subsequence) $\gamma_{k} \rightarrow \gamma$ in the varifold sense (and thus in Hausdorff distance). The limit varifold $\gamma$ is a connected, integral and stationary $1$-varifold.\\ \newline
\textbf{Claim $1$ : $\{ \gamma_{k} \}$ smoothly converges to $\gamma$ on a small neighborhood of each $x \in \gamma$. } \\
For each $k$, the intersection $B_{\delta}(x) \cap \gamma_{k}$ is a geodesic segment with zero curvature for sufficiently small $\delta$. Compactness argument gives that $\{ \gamma_{k} \}$ locally and smoothly converges to geodesic segments passing $x$.

Suppose there are more than one geodesic segments passing $x$ as a limit geodesic. Since $\gamma_{k}$ is simple, a segment of $\gamma_{k}$ converges to more than one piecewise curve touching at $x$. The geodesic curvature of each segment of piecewise curves touching $x$ blows up at $x$. It contradicts to the fact that the geodesic curvature of $\gamma$ is zero. Hence there is one geodesic segment passing $x$ in $B_{\delta}(x)$.

From above, $B_{\delta}(x) \cap \gamma$ consists of disjoint geodesic segments. By the Theorem (5) of Section $3$ in \cite{AA}, there are finite number of segments in $B_{\delta/2} (x) \cap \gamma$. By the compactness argument for each segment, we have a smooth and graphical convergence at each point $x \in \gamma$. \\ \newline
\textbf{Claim $2$ : If $M$ is an orientable surface, then the multiplicity of the convergence is $1$.} \\
By constancy theorem (Section $41$ of \cite {Si}), the density of $\gamma$ is a constant positive integer along $\gamma$. Suppose the convergence has a higher multiplicity $m>1$. Since $\gamma_{k}$ converges to $\gamma$ in a varifold sense, so does in Hausdorff distance sense hence $\gamma_{k}$ is within a $\delta$-normal neighborhood of $\gamma$ for small $\delta$ and sufficiently large $k$. Moreover, the $\delta$-normal neighborhood of $\gamma$ is diffeomorphic to $S^{1} \times (0,1)$ since $M$ is orientable. On a $\delta$-normal neighborhood $N_{\delta}(\alpha)$ of a small closed single segment $\alpha :=  B_{\delta}(x) \cap \gamma$, there is a set of functions $\{ u_{k}^{1} <u_{k}^{2} <...< u_{k}^{m}\} \in C^{\infty}(\alpha)$, $m>1$ such that
\begin{equation*}
   \gamma_{k} \cap N_{\delta}(\alpha) =\bigcup_{y \in \gamma} \{ Exp_{y}(\nu(y)u_{k}^{1}(y)), Exp_{y}(\nu(y)u_{k}^{2}(y)),..., Exp_{y}(\nu(y)u_{k}^{m}(y))\},
\end{equation*}
where $y \in B_{\delta}(x) \cap \gamma$ and $\nu(y)$ is a normal vector at $y$. Since $\gamma_{k}$ is locally $k$ disjoint segments for small neighborhood of each $x \in \gamma$ and $\gamma_{k}$ does not have any self-intersection point, graphs are separated to
\begin{equation*}
   \gamma_{k} = \bigcup_{i=1}^{m} \Big( \bigcup_{x \in \gamma} \{ Exp_{x}(\nu(x)u_{k}^{i}(x))\} \Big),
\end{equation*}
which are $k$ disjoint graphs. It contradicts to the fact that $\gamma_{k}$ is a single simple closed curve. Thus, the multiplicity of the convergence is $1$.
\\ \newline
\textbf{Claim 3 : If $M$ is a nonorientable surface, then the multiplicity of the convergence is at most $2$.} \\
We use the same notation as Claim 2 for a $\delta$-normal neighborhood and graphs. Suppose the convergence has a higher multiplicity $m>2$. In this case, the $\delta$-normal neighborhood of $\gamma$ is diffeomorphic to either a cylinder or a Möbius band.  The former case is the same as claim $2$ and the multiplicity is $1$. By the similar argument with Claim 2, the piece of curve containing $Exp_{x}(\nu(x)u_{k}^{1}(x))$ and $ Exp_{x}(\nu(x)u_{k}^{m}(x))$ is separated with other components of the curve due to the fact that the $N_{\delta}(\gamma)$ is diffeomorphic to a Möbius band and since $\gamma_{k}$ is locally $k$ disjoint segments in a small neighborhood. This contradicts to the fact that $\gamma_{k}$ is a single curve. Hence the multiplicity is at most $2$. \\ \newline
\textbf{Claim 4 : If $\gamma_{k} \cap \gamma = \varnothing$ eventually then $\gamma$ is stable, and $index(\gamma) \ge 1$ otherwise.} \\
This is due to the existence or nonexistence of strictly positive or negative eigenfunction and this argument is the same with the corresponding argument in \cite{Sh}.
\end{proof}
\end{thm}
\begin{cor}
On $(S^{2},g)$ with a bumpy metric, the set of all stationary integral varifolds $W_{L}$ in the sphere whose support is a simple closed, smooth geodesic and mass is less than or equal to $L$ is finite. 
\begin{proof}
Suppose $W_{L}$ is an infinite set, by the previous compactness theorem there exists an infinite sequence $\{w_{k} \}$ converging to $w$ in $W_{L}$. This induces a nontrivial Jacobi field of $w$ and contradicts to the fact that the metric is bumpy.
\end{proof}
\end{cor}

\end{document}